\documentclass[12pt,a4paper]{article}
\pdfoutput=1
\usepackage{amsmath}
\usepackage{amsfonts}
\usepackage{amssymb}
\usepackage{amsthm}
\usepackage{graphicx}

\newtheorem{theorem}{Theorem}[section]
\newtheorem{lemma}[theorem]{Lemma}
\newtheorem{corollary}[theorem]{Corollary}
\newtheorem{proposition}[theorem]{Proposition}
\newtheorem{definition}[theorem]{Definition}
\newtheorem{example}[theorem]{Example}

\begin{document}

\title{Continuous-Time Linear Optimization}

\author{Valeriano Antunes de Oliveira\thanks{This work was supported by the S\~ao Paulo Research Foundation (FAPESP) -- through the Center for Mathematical Sciences Applied to Industry -- CeMEAI/CEPID, under Grant 2013/07375-0.}}

\date{\small{S\~ao Paulo State University (Unesp), \\ Institute of Biosciences, Humanities and Exact Sciences, \\ S\~ao Jos\'e do Rio Preto, S\~ao Paulo, Brazil. \\ Email: valeriano.oliveira@unesp.br}}

\maketitle

\begin{abstract}
In this work, optimality conditions and classical results from duality theory are derived for continuous-time linear optimization problems with inequality constraints. The optimality conditions are given in the Karush-Kuhn-Tucker form. Weak and strong duality properties, as well as, the complementary slackness theorem are established. A result concerning the existence of solutions is also stated. \medskip

\end{abstract}

\section{Introduction}

This work is concerned with the continuous-time linear optimization problem posed as follows:
$$
\begin{array}{ll} 
\mbox{minimize} & F(z) = \displaystyle\int_{0}^{T} c(t)^\top z(t) dt \\ 
\mbox{subject to} & A(t)z(t) \leq b(t) ~ \mbox{a.e. in} ~ [0,T], \\
& z \in L^{\infty}([0,T];\mathbb{R}^n),
\end{array}
\eqno{\mathrm{(CLP)}}
$$
where $A : [0,T] \rightarrow \mathbb{R}^{m \times n}$, $b : [0,T] \rightarrow \mathbb{R}^{m}$ and $c : [0,T] \rightarrow \mathbb{R}^n$ have essentially bounded and measurable entries in $[0,T]$.

This class of problems was introduced in 1953 by Bellman \cite{bellman:1953} and since then the theory was considerably developed. Optimality conditions, as well as, duality results were obtained. At first, the matrix $A$ and vectors $b$ and $c$ were not allowed to vary with $t$, and the optimality conditions and duality results were proved to be valid under very strong assumptions. More general problems were then tackled and less restrictive assumptions on the problem data were made. In 1980, the article \cite{reiland:1980c} was published by Reiland generalizing previous results encountered in the literature until then. After the paper by Reiland, Zalmai published a series of articles on continuous-time nonlinear programming involving necessary and sufficient optimality conditions and duality theory. See \cite{zalmai:1985a}, for example. The assumptions made by Zalmai were less restrictive than those by Reiland. One of the main tools used by Zalmai was a generalization of the Gordan Transposition Theorem for the continuous-time context (see \cite{zalmai:1985d}). After the works by Zalmai, almost all the literature produced in this area made use of Gordan's Theorem. However, later, this theorem was discovered to be not valid (see \cite{arutyunov:2017}). Since then, the paper by Reiland became one of the main references on continuous-time linear programming. 

In \cite{reiland:1980c}, Reiland proved classical duality theorems for continuous-time linear programming problems under some constraint qualifications. He also presented an example showing that constraint qualifications are essential in establishing such results. So, in the continuous-time context, the linearity of the problem data is not itself a constraint qualification, as is the case on finite dimensions.

Here, we obtain optimality conditions and some classical results from duality theory for (CLP), namely, Karush-Kuhn-Tucker-type conditions, weak and strong duality properties, and the complementary slackness theorem. As pointed by Reiland, in the continuous-time context, even for linear problems, a constraint qualification is necessary. We propose a new constraint qualification that is less restrictive and simpler to be verified in comparison with the one used by Reiland. In fact, Reiland in \cite{reiland:1980c}, required (i) the validity of a Slater-type condition; (ii) that the kernel of a certain operator (between infinite-dimensional spaces) has finite dimensions and that its range is closed; and (iii) that the closure (in $L^\infty$) of the feasible directions cone coincides with the linearized feasible directions cone. The new constraint qualification proposed here stands only on a property of the cone generated by the rows of matrix $A(t)$. It is worth mentioning that our constraint qualification may be verified even when matrix $A(t)$ does not have full rank. However, in the obtained necessary optimality conditions, we may have a nonzero Lagrange multiplier related to a non-binding constraint. Indeed, we make use of the so-called $\beta$-active constraints ($\beta$ is a positive constant, see details in the next section). To the best of our knowledge, this notion appeared first in Osmolovski\u{\i} \cite{Osmolovski:1975} in the context of optimal control. Let us explain why it is necessary to work with $\beta$-active constraints. In finite dimensions, when one has $g(\bar{x}) < 0$, in which $g$ is continuous, then it holds $g(x) < 0$ for all $x$ in a certain neighborhood of $\bar{x}$. In the continuous-time context, if $g(\bar{z}(t),t) < 0$, then under the continuity of $g(\cdot,t)$, likewise $g(z(t),t) < 0$ for all $z(t)$ close enough to $\bar{z}(t)$. The point is that the neighborhood may vary with $t$. So, working with $\beta$-active constraints, that is, $g(\bar{z}(t),t) < -\beta$, $\beta > 0$, we obtain a neighborhood that depends only on $\bar{z} \in L^{\infty}([0,T];\mathbb{R}^n)$. In order to bypass the possibility of having nonzero Lagrange multipliers related to non-binding constraints, we also propose a full rank-type constraint qualification. Although stronger than the aforementioned new constraint qualification, the full rank condition may be useful when compared with Reiland's (i)-(iii) conditions mentioned above.

Roughly speaking, in Reiland \cite{reiland:1980c}, the results were obtained by specializing those given in Reiland \cite{reiland:1980b} for linear problems. We could follow this route, by specializing the Karush-Kuhn-Tucker-type necessary optimality conditions given in the recently published papers by Monte and de Oliveira \cite{monte:2019,monte:2020,monte:2021}. Notwithstanding, the approach employed here was direct. Moreover, the constraint qualifications used here are more adequate for linear programs. Compared to the full rank condition given in \cite{monte:2019}, in the full rank condition presented in this paper, only the $\beta$-active constraints are involved while in \cite{monte:2019} all constraints are taken into account. See more details in the next section. In \cite{monte:2020}, the optimality conditions are obtained under a Mangasarian-Fromovitz-type constraint qualification. Nevertheless, provided a theorem of the alternative is applied in the proof of the Karush-Kuhn-Tucker theorem, another regularity condition is assumed. This other regularity condition is specifically linked to the alternative theorem (see \cite{arutyunov:2017}). A constant rank-type constraint qualification is considered in \cite{monte:2021}. Although it is weaker than the full rank condition in \cite{monte:2019}, likewise, it involves all the constraints (not only the active or $\beta$-active ones). Further, the rank should be maintained not only in a neighborhood of the reference point $\bar{z} \in L^{\infty}([0,T];\mathbb{R}^n)$, but also with respect to $t$.

It should be mentioned that in Reiland \cite{reiland:1980c}, the continuous-time linear programming problem differs slightly from (CLP). The right hand side of the inequality constraint $A(t)z(t) \leq b(t)$ allows the presence of an integral term of the form $\int_t^T K(t,s)z(s)ds$, where $K(t,s)$ is a $m \times n$ matrix, and the feasible solutions should satisfy the non-negativity constraint $z(t) \geq 0$, $t \in [0,T]$.

Additionally, this paper brings a result concerning the existence of optimal solutions for (CLP). In fact, by assuming that the set of all feasible solutions is bounded, we show that (CLP) has an optimal solution. In Levinson \cite{levinson:1966}, it is proved a theorem on the existence of optimal solutions for a class of continuous-time linear programming problems. It is assumed that the entries of $A(t)$ are non-negative and piecewise continuous. The entries of $c(t)$ and $b(t)$ are also assumed to be piecewise continuous. Moreover, the existence theorem is proved under the hypothesis that there exist a piecewise continuous function $\lambda(t) \geq 0$, with $\| \lambda(t) \| \leq 1$, and a constant $b > 0$ such that $\| \lambda(t)^\top A(t) \| \geq b > 0$, $t \in [0,T]$. On the other hand, the problem considered by him is more general, similar to the one considered by Reiland. In the integral term, the entries of $K(t,s)$ are assumed to be piecewise continuous. The results by Levinson received important improvements through a perturbation approach given in Wu \cite{Wu:2013}. Wu treats the inequality constraints in the sense of a.e. in $[0, T]$ while Levinson assumes that they should be satisfied for all $t \in [0,T]$. A recent reference where an existence result is treated is Ghate \cite{Ghate:2020}. His result is similar to ours once it is only assumed that the feasible set is bounded, but the continuous-time problem is distinct: the coefficient matrix is not time-dependent, there is a non-negativity constraint, and the presence of an integral term similar to the one cited above but with a constant matrix.

Concerning the duality results, we obtain strong duality (with zero duality gap) and complementary slackness by assuming merely a full rank-type condition involving the $\beta$-active constraints. Among others, strong duality is also obtained by Levinson \cite{levinson:1966}, Reiland \cite{reiland:1980b}, Wu \cite{Wu:2013} and Zalmai \cite{zalmai:1985a}. In Levinson's paper, it is assumed that $A(t)$, $b(t)$ and $c(t)$ have continuous entries, there exists $\delta > 0$ such that either $a_{ij}(t) = 0$ or $a_{ij} \geq \delta$ for all $i,j$ and all $t \in [0,T]$, where $A(t) = (a_{ij})_{m \times n}$, $a_{ij}(t) \geq 0$, $t \in [0,T]$, and for each $t$ and $j$ there exists $i_j = i_j(t)$ such that $a_{i_jj}(t) \geq \delta$. Wu improves the result by assuming piecewise continuity instead of continuity and treating the inequality constraints in the sense of a.e. in $[0, T]$. In the strong duality theorem given in \cite{levinson:1966} and \cite{Wu:2013}, the existence of optimal solutions to both primal and dual problems is also guaranteed. In Reiland's result, in addition to the (i)-(iii) conditions mentioned above, assumptions similar to those by Levinson are necessary. On the other hand, the continuity hypotheses are unnecessary. In \cite{zalmai:1985a}, Zalmai proved the strong duality theorem under the assumption that the problem is stable, which, in turn, is equivalent to a certain perturbation function being finite as well as having nonempty subdifferential (in the sense of convex analysis) at the origin. Zalmai's results were developed in the more general setting of nonlinear problems (which, obviously, subsumes the linear case).

Considering that the problem is linear, and consequently convex, a possible approach to obtain duality results is through Rockafellar duality theory. See \cite{Rocka}. Similar results to those derived here can be achieved but under different assumptions. For example, in order to establish strong duality, while in \cite{Rocka} the assumption involves the subdifferential (in the sense of convex analysis) of a value function, here we have a full rank hypothesis on matrix $A(t)$. If, on the one hand, the assumptions in Rockafellar's book may be more general, on the other hand, ours are simpler to verify (from a practical point of view).

The work is organized in the following way. Next, some preliminaries are given. The main results concerning the optimality conditions are stated and proved in Section \ref{OptimCond} while those related to duality in Section \ref{dualidade}. Conclusion words are provided in the last section.

\section{Preliminaries}

In this section, we set the notation and give some key technical and preliminary results. Some important definitions are also stated.

The set of all feasible solutions of (CLP) will be denoted as $\Omega$:
$$
\Omega = \{ z \in L^{\infty}([0,T];\mathbb{R}^n) : A(t)z(t) \leq b(t) ~ \mbox{a.e. in} ~ [0,T] \}. 
$$

We say that $\bar{z} \in \Omega$ is an optimal solution of (CLP) if
$$
F(\bar{z}) \leq F(z) ~ \forall z \in \Omega.
$$

We assume throughout the paper that there exists a constant $K > 0$ such that, for almost every $t \in [0,T]$,
$$
\| A(t) \| \leq K, \quad \| b(t) \| \leq K, \quad \| c(t) \| \leq K.
$$

We will denote $I = \{ 1,2,\ldots,m \}$.

Let $\beta > 0$ be a small scalar. Given $\bar{z} \in \Omega$, we will denote by $I_\beta(t)$ the index set of the $\beta$-active constraints at instant $t$, that is,
$$
I_\beta(t) = \{ i \in I : -\beta \leq a_i(t)^\top \bar{z}(t) - b_i(t) \leq 0 \},
$$
where $a_i(t)^\top$ denotes the $i$-th row of matrix $A(t)$.

We will denote by $I_0(t)$ the index set of the active constraints at instant $t$, that is,
$$
I_0(t) = \{ i \in I : a_i(t)^\top \bar{z}(t) - b_i(t)= 0 \}.
$$

Given an indices subset $\mathcal{I} \subset I$ and $t \in [0,T]$, we will denote by $A^{\mathcal{I}}(t)$ the matrix obtained from $A(t)$ by removing the rows whose indices do not belong to $\mathcal{I}$. Analogously for any vector or matrix to be defined throughout the paper. The cardinality of $\mathcal{I}$ will be denoted by $|\mathcal{I}|$. 

The singular values of $A^{\mathcal{I}}(t)$ will be denoted by $\sigma_{i}(A^{\mathcal{I}}(t))$, $i = 1, \ldots, q(t)$, where $q(t) = \min \{ |\mathcal{I}|,n \}$. Let us remember that the singular values are non-increasingly ordered, that is, 
$$
\sigma_{1}(A^{\mathcal{I}}(t)) \geq \sigma_{2}(A^{\mathcal{I}}(t)) \geq \cdots \geq \sigma_{q(t)}(A^{\mathcal{I}}(t)).
$$
The smallest positive singular value of $A^{\mathcal{I}}(t)$ will be denoted by $\sigma_{\min}(A^{\mathcal{I}}(t))$. It is known that $\sigma_{\min}(A^{\mathcal{I}}(t))$ coincides with the $p(t)$-th singular value of $A^{\mathcal{I}}(t)$, where $p(t) = \mathrm{rank}(A^{\mathcal{I}}(t))$.

The cone generated by the rows of $A^{\mathcal{I}}(t)$ will be denoted by $\mathrm{cone}(A^{\mathcal{I}}(t))$.


\begin{definition} \label{beta_CR}
Let $\beta > 0$. We say that the Regularity Condition ($\beta$-RC) is satisfied at $\bar{z} \in \Omega$ if there exist an indices subset $I_{\sharp}(t) \subset I_\beta(t)$ and a positive constant $\bar{K}$ such that
$$
\mathrm{cone}(A^{I_{\sharp}(t)}(t)) = \mathrm{cone}(A^{I_\beta(t)}(t)) ~ \mbox{a.e. in} ~ [0,T]
$$
and
$$
\mathrm{det}(A^{I_{\sharp}(t)}(t)A^{I_{\sharp}(t)}(t)^\top) \geq \bar{K} ~ \mbox{a.e. in} ~ [0,T].
$$
\end{definition}

\begin{definition} \label{beta_FR}
Let $\beta > 0$. We say that the Full Rank Condition ($\beta$-FR) is satisfied at $\bar{z} \in \Omega$ if there exists a positive constant $\hat{K}$ such that
$$
\mathrm{det}(A^{I_{\beta}(t)}(t)A^{I_{\beta}(t)}(t)^\top) \geq \hat{K} ~ \mbox{a.e. in} ~ [0,T].
$$
\end{definition}

It is clear that $\beta$-FR implies in $\beta$-CR with $I_\sharp = I_\beta$.

In Monte and de Oliveira \cite{monte:2019}, a distinct full rank condition is proposed in the context of nonlinear continuous-time programming. In the case of (CLP), it is satisfied at a feasible solution $\bar{z}$ if there exists $K_{\bar{z}} > 0$ such that such that $\det(\Gamma(t)\Gamma(t)^\top) \geq K_{\bar{z}}$ a.e. in $[0,T]$ where
$
\Gamma(t) = \left[
\begin{array}{cc}
A(t) & D(t)
\end{array}
\right],
$ 
$D(t) = \mathrm{diag}(\{ -2 \sqrt{a_i(t)^\top \bar{z}(t) - b_i(t)} \}_{i \in I})$. As pointed in Remark 1 in \cite{monte:2019}, $\Gamma(t)$ has full rank if, and only if, $A^{I_0(t)}(t)$ has full rank. Note that $A^{I_0(t)}(t)$ having full rank is not the same as $\det(A^{I_0(t)}(t)A^{I_0(t)}(t)^\top)$ being bounded below away from zero. On the other hand, if the full rank condition is valid, then $\det(A^{I_0(t)}(t)A^{I_0(t)}(t)^\top)$ is bounded below away from zero. The converse is not true. Moreover, after an application of Weyl Theorem (see Horn and Johnson \cite{Horn:2013}, for example), it is easy to see that the full rank condition given in \cite{monte:2019} is valid when $\det(A(t)A(t)^\top) \geq K_{\bar{z}}$ a.e. in $[0,T]$. Observe that all constraints but not only the $\beta$-active are involved in this last condition.

Before we go on, let us remember some useful properties of singular values of a matrix. Given a matrix $M \in \mathbb{R}^{m \times n}$, if $M = U \Sigma V^T$ is its singular values decomposition, the matrix $M^+ \in \mathbb{R}^{n \times m}$ defined as $M^+ = V \Sigma^+ U^T$ is called the Moore-Penrose pseudo-inverse of $M$, where 
\[
\Sigma^+ = \left[ \begin{array}{cc} E & 0 \\ 0 & 0 \end{array} \right] \in \mathbb{R}^{n \times m}, ~ E = \mathrm{diag}(\sigma_i(M)^{-1})_{i=1}^{p}, ~ p = \mathrm{rank}(M),
\]
and $\sigma_1(M) \geq \sigma_2(M) \geq \cdots \geq \sigma_p(M) > 0$ are the positive singular values of $M$. The following important properties are valid.
\begin{enumerate}
\item[(i)] $MM^+M = M$ and $M^+MM^+ = M^+$;
\item[(ii)] $(MM^+)^\top = MM^+$ and $(M^+M)^\top = M^+M$;
\item[(iii)] $(M^+)^\top = (M^\top)^+$ and $(M^+)^+ = M$;
\item[(iv)] $\| M \| = \sigma_1(M)$ and $\| M^+ \| = \sigma_p(M)^{-1}$;
\item[(v)] The nonzero singular values of $M$ coincide with the square roots of the nonzero eigenvalues of $MM^\top$ or $M^\top M$;
\item[(vi)] If $n=m$ and $M$ is symmetric, the singular values of $M$ coincide with the absolute values of the eigenvalues of $M$;
\item[(vii)] If $n=m$, then $|\det(M)| = \prod_{i=1}^{m} \sigma_i(M)$. 
\end{enumerate}
For more properties of the pseudo-inverse matrix and the singular value decomposition, the reader is referred to Horn and Johnson \cite{Horn:2013}, Noble and Daniel \cite{Noble:1977} or Trefethen and Bau \cite{Trefethen:1997}.

Following, we have some technical lemmas which will be useful later in the paper.

\begin{lemma} \label{Lema1}
The full rank condition $\beta$-FR is valid if, and only if, the matrix $A^{I_\beta(t)}(t)$ has full rank and there exists a positive constant $C$ such that \linebreak $\sigma_{\min}(A^{I_\beta(t)}(t)) \geq C$ a.e. in $[0,T]$.
\end{lemma}

\begin{proof}
Let us assume that $\beta$-FR is valid at $\bar{z} \in \Omega$. It is obvious that $A^{I_\beta(t)}(t)$ has full rank a.e. in $[0,T]$. We have, almost everywhere in $[0,T]$, that
$$
\sigma_i(A^{I_\beta(t)}(t)) \leq \| A^{I_\beta(t)}(t) \| \leq K ~ \forall i,
$$
$$
\sigma_i(A^{I_\beta(t)}(t)A^{I_\beta(t)}(t)^\top) = [\sigma_i(A^{I_\beta(t)}(t))]^2 ~ \forall i
$$
and
$$
\prod_{i = 1}^{\vert I_\beta(t) \vert} \sigma_i(A^{I_\beta(t)}(t)A^{I_\beta(t)}(t)^\top) = \det(A^{I_\beta(t)}(t)A^{I_\beta(t)}(t)^\top) \geq \hat{K}.
$$
Thence,
\begin{eqnarray*}
\sigma_{\vert I_\beta(t) \vert}(A^{I_\beta(t)}(t)A^{I_\beta(t)}(t)^\top) & \geq & \frac{\hat{K}}{\prod_{i = 1}^{\vert I_\beta(t) \vert -1} \sigma_i(A^{I_\beta(t)}(t)A^{I_\beta(t)}(t)^\top)} \\ 
& \geq & \frac{\hat{K}}{(K^2)^{\vert I_\beta(t) \vert -1}} ~ \mbox{a.e. in} ~ [0,T].
\end{eqnarray*}
Therefore,
\begin{eqnarray*}
\sigma_{\min}(A^{I_\beta(t)}(t)) & = & \sqrt{\sigma_{\min}(A^{I_\beta(t)}(t)A^{I_\beta(t)}(t)^\top)} \\
& = & \sqrt{\sigma_{\vert I_\beta(t) \vert}(A^{I_\beta(t)}(t)A^{I_\beta(t)}(t)^\top)} \\
& \geq & \sqrt{\frac{\hat{K}}{(K^2)^{\vert I_\beta(t) \vert -1}}}  ~ \mbox{a.e. in} ~ [0,T].
\end{eqnarray*}
By taking
$$
C := \left\{
\begin{array}{ll}
\sqrt{\frac{\hat{K}}{(K^2)^{m -1}}} & \mbox{if} ~ K^2 \geq 1, \\
\sqrt{\hat{K}} & \mbox{if} ~ 0 < K^2 < 1, 
\end{array}
\right.
$$
it follows that $\sigma_{\min}(A^{I_\beta(t)}(t)) \geq C$ a.e. in $[0,T]$.

Conversely, assume that $A^{I_\beta(t)}(t)$ has full rank and there exists a positive constant $C$ such that $\sigma_{\min}(A^{I_\beta(t)}(t)) \geq C$ a.e. in $[0,T]$. We have that
$$
\sigma_{i}(A^{I_\beta(t)}(t)) \geq \sigma_{\min}(A^{I_\beta(t)}(t)) \geq C ~ \mbox{a.e. in} ~ [0,T], ~ \forall i.
$$
Then,
\begin{eqnarray*}
\det(A^{I_\beta(t)}(t)A^{I_\beta(t)}(t)^\top) & = & \prod_{i = 1}^{\vert I_\beta(t) \vert} \sigma_i(A^{I_\beta(t)}(t)A^{I_\beta(t)}(t)^\top) \\
& = & \prod_{i = 1}^{\vert I_\beta(t) \vert} [\sigma_i(A^{I_\beta(t)}(t))]^2 \geq C^{\vert I_{\beta(t)} \vert} ~ \mbox{a.e. in} ~ [0,T].
\end{eqnarray*}
Defining
$$
\hat{K} := \left\{
\begin{array}{ll}
C & \mbox{if} ~ C \geq 1, \\
C^m & \mbox{if} ~ 0 < C < 1,
\end{array}
\right.
$$
we get $\det(A^{I_\beta(t)}(t)A^{I_\beta(t)}(t)^\top)  \geq \hat{K}$ a.e. in $[0,T]$.
\end{proof}

\begin{lemma} \label{Lema2}
If the regularity condition $\beta$-RC is valid, then there exists a positive constant $C$ such that $\sigma_{\min}(A^{I_\beta(t)}(t)) \geq C$ a.e. in $[0,T]$.
\end{lemma}

\begin{proof}
Let us denote $I_\flat(t) = I_\beta(t) \setminus I_\sharp(t)$ a.e. in $[0,T]$. By permuting the rows of $A^{I_\beta(t)}(t)$ in such a way that the rows whose indices belong to $I_\sharp(t)$ come first, we can write
$$
A^{I_\beta(t)}(t) = \left[
\begin{array}{c}
A^{I_\sharp(t)}(t) \\ A^{I_\flat(t)}(t)
\end{array}
\right] ~ \mbox{a.e. in} ~ [0,T].
$$
Then, 
$$
A^{I_\beta(t)}(t)^\top A^{I_\beta(t)}(t) = A^{I_\sharp(t)}(t)^\top A^{I_\sharp(t)}(t) + A^{I_\flat(t)}(t)^\top A^{I_\flat(t)}(t) ~ \mbox{a.e. in} ~ [0,T].
$$
We will apply Weyl Theorem (see Theorem 4.3.1 in Horn and Johnson \cite{Horn:2013}, for instance) in the equality above. Let us denote the eigenvalues of $A^{I_\beta(t)}(t)^\top A^{I_\beta(t)}(t)$ by $\lambda_i(A^{I_\beta(t)}(t)^\top A^{I_\beta(t)}(t))$ a.e. in $[0,T]$. Similarly, we use the same notation for $A^{I_\sharp(t)}(t)^\top A^{I_\sharp(t)}(t)$ and  $A^{I_\flat(t)}(t)^\top A^{I_\flat(t)}(t)$. Let us remember that, in Theorem 4.3.1 of \cite{Horn:2013}, the non-decreasing ordering of the eigenvalues is considered. If $r(t) := |I_\sharp(t)| = \mathrm{rank}(A^{I_\sharp(t)}(t)) = \mathrm{rank}(A^{I_\beta(t)}(t)) = \mathrm{rank}(A^{I_\beta(t)}(t)^\top A^{I_\beta(t)}(t))$ a.e. in $[0,T]$, by taking $i=j=n-r(t)+1$ in Weyl Theorem, we obtain
\begin{multline*}
\lambda_{n-r(t)+1}(A^{I_\beta(t)}(t)^\top A^{I_\beta(t)}(t))  \\
\geq \lambda_1(A^{I_\flat(t)}(t)^\top A^{I_\flat(t)}(t)) + \lambda_{n-r(t)+1}(A^{I_\sharp(t)}(t)^\top A^{I_\sharp(t)}(t)) ~ \mbox{a.e. in} ~ [0,T]. 
\end{multline*}
Provided $A^{I_\flat(t)}(t)^\top A^{I_\flat(t)}(t)$ is positive semi-definite, $\lambda_1(A^{I_\flat(t)}(t)^\top A^{I_\flat(t)}(t)) \geq 0$ a.e. in $[0,T]$, so that
\begin{multline*}
\lambda_1(A^{I_\flat(t)}(t)^\top A^{I_\flat(t)}(t)) + \lambda_{n-r(t)+1}(A^{I_\sharp(t)}(t)^\top A^{I_\sharp(t)}(t)) \\
\geq \lambda_{n-r(t)+1}(A^{I_\sharp(t)}(t)^\top A^{I_\sharp(t)}(t)) ~ \mbox{a.e. in} ~ [0,T].
\end{multline*}
Hence,
$$
\lambda_{n-r(t)+1}(A^{I_\beta(t)}(t)^\top A^{I_\beta(t)}(t)) \geq \lambda_{n-r(t)+1}(A^{I_\sharp(t)}(t)^\top A^{I_\sharp(t)}(t)) ~ \mbox{a.e. in} ~ [0,T]. 
$$
We know that the nonzero singular values of $A^{I_\beta(t)}(t)$ coincide with the square roots of the nonzero eigenvalues of $A^{I_\beta(t)}(t)^\top A^{I_\beta(t)}(t)$ a.e. in $[0,T]$. Provided the singular values are decreasingly ordered and the non-decreasing ordering of the eigenvalues is considered in Theorem 4.3.1, we have 
$$
\sigma_{\min}(A^{I_\beta(t)}(t))^2 = \sigma_{r(t)}(A^{I_\beta(t)}(t))^2 = \lambda_{n-r(t)+1}(A^{I_\beta(t)}(t)^\top A^{I_\beta(t)}(t)) ~ \mbox{a.e. in} ~ [0,T].
$$
Analogously,
$$
\sigma_{\min}(A^{I_\sharp(t)}(t))^2 = \sigma_{r(t)}(A^{I_\sharp(t)}(t))^2 = \lambda_{n-r(t)+1}(A^{I_\sharp(t)}(t)^\top A^{I_\sharp(t)}(t)) ~ \mbox{a.e. in} ~ [0,T].
$$
Hence, $\sigma_{\min}(A^{I_\beta(t)}(t))^2 \geq \sigma_{\min}(A^{I_\sharp(t)}(t))^2$ a.e. in $[0,T]$.
Proceeding as in the first part of the proof of Lemma \ref{Lema1}, we get that there exists a positive constant $C$ such that $\sigma_{\min}(A^{I_\sharp(t)}(t)) \geq C$ a.e. in $[0,T]$, from where the result follows. 
\end{proof}

Note that the cone condition $\mathrm{cone}(A^{I_{\sharp}(t)}(t)) = \mathrm{cone}(A^{I_\beta(t)}(t))$ a.e. in $[0,T]$ in Definition \ref{beta_CR} does not play any rule in the proof above.

\begin{lemma} \label{Lema3}
If the regularity condition $\beta$-RC is valid, then there exists a positive constant $\tilde{K}$ such that $\| [A^{I_\sharp(t)}(t) A^{I_\sharp(t)}(t)^\top]^{-1} \| \leq \tilde{K}$ a.e. in $[0,T]$.
\end{lemma}

\begin{proof}
As already mentioned above, proceeding as in the first part of the proof of Lemma \ref{Lema1}, we get that there exists a positive constant $C$ such that $\sigma_{\min}(A^{I_\sharp(t)}(t)) \geq C$ a.e. in $[0,T]$. We have that
\begin{eqnarray*}
\| [A^{I_\sharp(t)}(t) A^{I_\sharp(t)}(t)^\top]^{-1} \| & = & \frac{1}{\sigma_{\min}(A^{I_\sharp(t)}(t) A^{I_\sharp(t)}(t)^\top)} \\
& = & \frac{1}{\sigma_{\min}(A^{I_\sharp(t)}(t))^2} \leq \frac{1}{C^2} ~ \mbox{a.e. in} ~ [0,T].
\end{eqnarray*}

\end{proof}



Below we present an auxiliary result, which will be used in the proof of some results in Section \ref{OptimCond}. Let $\phi : \mathbb{R}^n \times [0,T] \rightarrow \mathbb{R}$. Assume that (i) $\phi(\cdot,t)$ is continuously differentiable throughout $[0,T]$, (ii) $\phi(z,\cdot)$ is measurable for each $z$ and (iii) there exists $K_\phi > 0$ such that $\| \nabla_z \phi(\bar{z}(t),t) \| \leq K_\phi$ a.e. in $[0,T]$.

\begin{proposition} \label{prob-sem-restr}
Let $\bar{z} \in L^{\infty}([0,T];\mathbb{R}^n)$ be an optimal solution to 
$$
\begin{array}{ll} 
\mbox{minimize} & F(z) = \displaystyle\int_{0}^{T} \phi(z(t),t) dt \\ 
\mbox{subject to} & z \in L^{\infty}([0,T];\mathbb{R}^n).
\end{array}
$$
Then, 
$$
\nabla_z \phi(\bar{z}(t),t) = 0 ~ \mbox{a.e. in} ~ [0,T].
$$
\end{proposition} 

The proof can be found in Monte and de Oliveira \cite{monte:2019}. Actually, this is a particular case of Proposition 3.2 in \cite{monte:2019}.

We finish this section with an existence result.
\begin{proposition} \label{Existencia}
Assume that $\Omega \neq \emptyset$ and $\Omega$ is bounded. Then (CLP) has an optimal solution.
\end{proposition}

The result follows after applying the classical Weierstrass Theorem. One should sees $L^\infty([0,T];\mathbb{R}^n)$ as the dual space of $L^1([0,T];\mathbb{R}^n)$ and uses the weak$^\ast$ topology.

\section{Optimality Conditions} \label{OptimCond}

This section is devoted to establishing the optimality conditions for (CLP). The results stated next give the necessary conditions in the Karush-Kuhn-Tucker form. Such conditions are shown to be also sufficient for optimality, which is natural given the linearity of the problem. An illustrative example is given at the end of the section.

\begin{theorem} \label{CondNecPRL}
Let $\bar{z} \in \Omega$ be an optimal solution of (CLP). Assume that $\beta$-RC is satisfied at $\bar{z}$ for some $\beta > 0$. Then, there exists $u \in L^{\infty}([0,T];\mathbb{R}^m)$ such that
\begin{eqnarray}
&& c(t) + \sum_{i \in I_\sharp(t)} u_i(t)a_i(t) = 0 ~ \mbox{a.e. in} ~ [0,T], \label{KKT_PL} \\
&& u_i(t) \geq 0  ~ \mbox{a.e. in} ~ [0,T], ~ i \in I. \label{NNeg_PL}
\end{eqnarray}
\end{theorem}

\begin{proof}
For almost every $t \in [0,T]$ and all $i \in I$, we denote
$$
\bar{g}_i(t) = a_i(t)^\top \bar{z}(t) - b_i(t) \qquad \mbox{and} \qquad
\delta_i^\beta(t) = \left\{
\begin{array}{l}
1 ~ \mbox{if} ~ i \in I_\beta(t), \\ 0 ~ \mbox{otherwise}.
\end{array}
\right.
$$
Note that $\bar{g}_i(t) \leq 0$ a.e. in $[0,T]$. Assume that system below has a solution $\gamma \in L^{\infty}([0,T];\mathbb{R}^n)$.
\begin{equation} \label{sistema2}
\left\{ \begin{array}{l}
\int_0^T c(t)^\top \gamma(t) dt < 0, \\ \delta_i^\beta(t) a_i(t)^\top \gamma(t) \leq 0, ~ i \in I, ~ \mbox{a.e. in} ~ [0,T].
\end{array} \right.
\end{equation}
Let $K_1$ and $K_2$ be positive scalars such that $|\gamma(t)| \leq K_1$ and $\| a_i(t) \| \leq K_2$, $i \in I$, a.e. in $[0,T]$. Let $0< \tau < \beta K_1^{-1}K_2^{-1}$. We have that
\[
F(\bar{z}+\tau\gamma) - F(\bar{z}) = \tau \int_0^T c(t)^\top \gamma(t) dt < 0,
\] 
so that, $F(\bar{z}+\tau\gamma) < F(\bar{z})$. Further, for almost every $t \in [0,T]$, by considering $i \in I \setminus I_\beta(t)$, we have that 
\begin{eqnarray*}
a_i(t)^\top [\bar{z}(t) + \tau \gamma(t)] - b_i(t) & = & [a_i(t)^\top \bar{z}(t) - b_i(t)] + \tau a_i(t)^\top \gamma(t) \\
& < & - \beta + \tau \| a_i(t) \| \| \gamma(t) \| \leq - \beta + \tau K_2K_1 < 0.
\end{eqnarray*}
Then $a_i(t)^\top [\bar{z}(t) + \tau \gamma(t)] < b_i(t)$ a.e. in $[0,T]$. If $i \in I_\beta(t)$, it follows from the second inequality in (\ref{sistema2}) that $a_i(t)^\top \gamma(t) \leq 0$, from where we see that
\begin{eqnarray*}
a_i(t)^\top [\bar{z}(t) + \tau \gamma(t)] - b_i(t) & = & [a_i(t)^\top \bar{z}(t) - b_i(t)] + \tau a_i(t)^\top \gamma(t) \\
& = & \bar{g}_i(t) + \tau a_i(t)^\top \gamma(t) \leq 0.
\end{eqnarray*}
Thence, $a_i(t)^\top [\bar{z}(t) + \tau \gamma(t)] \leq b_i(t)$ a.e. in $[0,T]$. 

Summarizing, we obtain $\bar{z} + \tau \gamma \in \Omega$ with $F(\bar{z}+\tau\gamma) < F(\bar{z})$ for all $\tau > 0$ small enough, which contradicts the optimality of $\bar{z}$ in (CLP). Therefore, system (\ref{sistema2}) does not have any solution.

For almost every $t \in [0,T]$, we denote 
$$
\Delta^\beta(t) = \mathrm{diag}(\delta_1^\beta(t),\ldots,\delta_m^\beta(t)) \qquad \mbox{and} \qquad M(t) = \Delta^\beta(t) A(t).
$$
Provided system (\ref{sistema2}) does not have any solution, it is easy to see that $\bar{\gamma}=0$ is an optimal solution to the continuous-time optimization problem below. 
\begin{equation} \label{problema_auxiliar1}
\begin{array}{ll} 
\mbox{Minimize} & F(\gamma) = \displaystyle\int_{0}^{T} c(t)^\top \gamma(t) dt \\ 
\mbox{subject to} & M(t) \gamma(t) \leq 0 ~ \mbox{a.e. in} ~ [0,T], \\
& \gamma \in L^{\infty}([0,T];\mathbb{R}^n).
\end{array}
\end{equation}
Let us consider the auxiliary unconstrained continuous-time optimization problem posed as follows:
\begin{equation} \label{problema_auxiliar2}
\begin{array}{ll} 
\mbox{minimize} & \Phi(\gamma) = \displaystyle\int_{0}^{T} c(t)^\top \left[ \gamma(t) - M(t)^+M(t)\gamma(t) \right] dt \\ 
\mbox{subject to} & \gamma \in L^{\infty}([0,T];\mathbb{R}^n).
\end{array}
\end{equation}
Let $\gamma \in L^{\infty}([0,T];\mathbb{R}^n)$ and set 
$$
\hat{\gamma}(t) = \gamma(t) - M(t)^+M(t) \gamma(t) ~ \mbox{a.e. in} ~ [0,T].
$$
We know that the pseudo-inverse matrix has the properties that $|M(t)^+M(t)| = 1$ and $M(t)M(t)^+M(t)=M(t)$. So, it is clear that $\hat{\gamma} \in L^{\infty}([0,T];\mathbb{R}^n)$ and
$$
M(t) \hat{\gamma}(t) = 0 ~ \mbox{a.e. in} ~ [0,T].
$$
We see that $\hat{\gamma}$ is a feasible solution to problem (\ref{problema_auxiliar1}). Then, $F(\hat{\gamma}) \geq F(\bar{\gamma})$, so that
$$
\Phi(\gamma) = F(\hat{\gamma}) \geq F(\bar{\gamma}) = \Phi(\bar{\gamma}).
$$
Thus, $\bar{\gamma}$ is an optimal solution of problem (\ref{problema_auxiliar2}). It follows from Proposition \ref{prob-sem-restr} that
$$
c(t) - M(t)^\top (M(t)^+)^\top c(t) = 0 ~ \mbox{a.e. in} ~ [0,T].
$$
Defining
$$
\tilde{u}(t) = -(M(t)^+)^\top c(t) ~ \mbox{a.e. in} ~ [0,T],
$$
we have that
\begin{equation} \label{eqt}
c(t) + A(t)^\top \Delta^\beta(t) \tilde{u}(t) = 0 ~ \mbox{a.e. in} ~ [0,T].
\end{equation}
Moreover, $\tilde{u} \in L^{\infty}([0,T];\mathbb{R}^m)$. Indeed, if $K$ is such that $|c(t)| \leq K$ a.e. in $[0,T]$,
\begin{eqnarray}
\| \tilde{u}(t) \| &\leq& \| (M(t)^+)^\top \| \|c(t)\| = \| M(t)^+ \| \| c(t) \| \nonumber \\ 
&=& \sigma_{\min}(M(t))^{-1} |c(t)| \nonumber \\ 
& = & \sigma_{\min}(A^{I_\beta(t)}(t))^{-1} |c(t)| \leq C^{-1}K ~ \mbox{a.e. in} ~ [0,T], \label{eqt3}
\end{eqnarray}
where we used Lemma \ref{Lema2} in the last inequality.

Let us denote $I_\flat(t) = I_\beta(t) \setminus I_\sharp(t)$ a.e. in $[0,T]$. From the assumption that $\beta$-RC is satisfied at $\bar{z}$, it follows that there exists a matrix $\Lambda(t) = [\lambda_{ij}(t)] \in \mathbb{R}^{|I_\flat(t)| \times |I_\sharp(t)|}$ with $\lambda_{ij}(t) \geq 0$ for all $i,j$ such that 
\begin{equation} \label{eqt2}
A^{I_\flat(t)}(t) = \Lambda(t) A^{I_\sharp(t)}(t) ~ \mbox{a.e. in} ~ [0,T].
\end{equation}
It follows also that 
$A^{I_\sharp(t)}(t)A^{I_\sharp(t)}(t)^\top$ is invertible a.e. in $[0,T]$. We have that
$$
A^{I_\flat(t)}(t) A^{I_\sharp(t)}(t)^\top = \Lambda(t) A^{I_\sharp(t)}(t) A^{I_\sharp(t)}(t)^\top,
$$
which implies in
$$
\Lambda(t) = A^{I_\flat(t)}(t) A^{I_\sharp(t)}(t)^\top [A^{I_\sharp(t)}(t) A^{I_\sharp(t)}(t)^\top]^{-1} ~ \mbox{a.e. in} ~ [0,T].
$$
From Lemma \ref{Lema3}, it follows that there exist $\tilde{K} > 0$ such that 
$$
\| [A^{I_\sharp(t)}(t) A^{I_\sharp(t)}(t)^\top]^{-1} \| \leq \tilde{K} ~ \mbox{a.e. in} ~ [0,T]. 
$$
Thence,
\begin{equation} \label{eqt1}
\| \Lambda(t) \| \leq \tilde{K} \| A(t) \|^2 \leq \tilde{K} K^2 ~ \mbox{a.e. in} ~ [0,T].
\end{equation}

Set, for almost every $t \in [0,T]$, $\bar{u}(t) = \tilde{u}^{I_\sharp(t)}(t) + \Lambda(t)^\top \tilde{u}^{I_\flat(t)}(t) \in \mathbb{R}^{|I_\sharp(t)|}$ and $u(t) \in \mathbb{R}^m$ in which
$$
u_i(t) = \left\{ \begin{array}{ll}
\bar{u}_i(t), & \mbox{for} ~ i \in I_\sharp(t), \\ 0, & \mbox{otherwise}.
\end{array} \right.
$$
We have from (\ref{eqt}) and (\ref{eqt2}) that
\begin{eqnarray}
c(t) + A(t)^\top \Delta^\beta(t) u(t) &=& c(t) + A^{I_\sharp(t)}(t)^\top \bar{u}(t) \nonumber \\
&=& c(t) +A^{I_\sharp(t)}(t)^\top [ \tilde{u}^{I_\sharp(t)}(t) + \Lambda(t)^\top \tilde{u}^{I_\flat(t)}(t) ] \nonumber \\
&=& c(t) +A^{I_\sharp(t)}(t)^\top \tilde{u}^{I_\sharp(t)}(t) + A^{I_\flat(t)}(t)^\top \tilde{u}^{I_\flat(t)}(t) \nonumber \\
&=&c(t) + A(t)^\top \Delta^\beta(t) \tilde{u}(t) \nonumber \\
&=& 0 ~ \mbox{a.e. in} ~ [0,T] \label{eqt6}
\end{eqnarray}
and, from (\ref{eqt3}) and (\ref{eqt1}),
\begin{eqnarray*}
\| u(t) \| = \| \bar{u}(t) \| & \leq & \| \tilde{u}^{I_\sharp(t)}(t) \| + \| \Lambda(t) \| \| \tilde{u}^{I_\flat(t)}(t) \| \\
& \leq & C^{-1}K + \tilde{K} K^2 C^{-1}K ~ \mbox{a.e. in} ~ [0,T].
\end{eqnarray*}
The last inequality shows that $u \in L^{\infty}([0,T];\mathbb{R}^m)$. 

Condition (\ref{KKT_PL}) follows from (\ref{eqt6}) and the definition of $u$.

We now show that $u_i(t) = \bar{u}_i(t) \geq 0$ for $i \in I_\sharp(t)$ a.e. in $[0,T]$. Suppose that there exists $S \subset [0,T]$ with positive measure such that $u_{i_0}(t) < 0$ for some $i_0 \in I^\sharp(t)$, for all $t \in S$. Define, almost everywhere in $[0,T]$, $d(t) \in \mathbb{R}^{m}$ as $d_i(t) = 0$ for $i \in I_\sharp(t) \cup (I \setminus I_\beta(t))$, $i \neq i_0$,
$$
d_{i_0}(t) = \left\{ \begin{array}{rl}
-1, & \mbox{if} ~ t \in S, \\ 0, & \mbox{otherwise},
\end{array} \right.
$$
and
$$
d_i(t) = \sum_{j \in I_\sharp(t)} \lambda_{ij}(t)d_j(t), ~ i \in I_\flat(t), ~\mbox{a.e. in} ~ [0,T].
$$
Note that $d^{I_\flat(t)}(t) = \Lambda(t) d^{I_\sharp(t)}(t)$. It is easy to see from (\ref{eqt1}) that $d \in L^{\infty}([0,T];\mathbb{R}^m)$. Take 
$$
\tilde{\gamma}(t) = A^{I_\sharp(t)}(t)^\top [A^{I_\sharp(t)}(t) A^{I_\sharp(t)}(t)^\top]^{-1} d^{I_\sharp(t)}(t) ~\mbox{a.e. in} ~ [0,T].
$$
It is clear that $\tilde{\gamma} \in L^{\infty}([0,T];\mathbb{R}^n)$. We have also that
\begin{equation} \label{eqt4}
A^{I_\sharp(t)}(t) \tilde{\gamma}(t) = d^{I_\sharp(t)}(t) \leq 0,
\end{equation}
and
\begin{equation} \label{eqt5}
A^{I_\flat(t)}(t) \tilde{\gamma}(t) = \Lambda(t) A^{I_\sharp(t)}(t) \tilde{\gamma}(t) = \Lambda(t) d^{I_\sharp(t)}(t) \leq 0
\end{equation}
a.e. in $[0,T]$, where in the last inequality we used that fact that the entries of $\Lambda(t)$ are non-negative and $d^{I_\sharp(t)}(t) \leq 0$. Therefore, $M(t) \tilde{\gamma}(t) \leq 0$ a.e. in $[0,T]$, that is, $\tilde{\gamma}$ is a feasible solution to (\ref{problema_auxiliar1}). It follows that $F(\tilde{\gamma}) \geq F(\bar{\gamma}) = 0$. On the other hand, by using (\ref{eqt6}), (\ref{eqt4}) and (\ref{eqt5}),
\begin{eqnarray*}
F(\tilde{\gamma}) &=& \int_0^T c(t)^\top \tilde{\gamma}(t) dt = -\int_0^T u(t)^\top \Delta^\beta(t) A(t) \tilde{\gamma}(t) dt =  -\int_0^T u(t)^\top d(t) dt \\
&=& -\int_0^T \bar{u}(t)^\top d^{I_\sharp(t)}(t) dt = -\int_0^T u_{i_0}(t)d_{i_0}(t) dt = \int_S u_{i_0}(t) dt < 0. 
\end{eqnarray*}
So, we come to a contradiction. Therefore, condition (\ref{NNeg_PL}) is valid.
\end{proof}

In Definition \ref{beta_CR}, note that though $I_\sharp(t) \subset I_\beta(t)$ it is not possible to assure that $I_\sharp(t) \subset I_0(t)$. So, in the optimality conditions obtained in Theorem \ref{CondNecPRL}, we may have $u_i(t) \neq 0$ for $i \notin I_0(t)$, that is, a multiplier related to a non-biding (though $\beta$-active) constraint can be nonzero. As already mentioned, when $\beta$-FR is valid, $\beta$-CR is fulfilled with $I_\sharp(t) = I_\beta(t)$. In this case, the inconvenience pointed above disappears.

\begin{theorem} \label{CondNecPRL2}
Let $\bar{z} \in \Omega$ be an optimal solution of (CLP). Assume that $\beta$-FR is satisfied at $\bar{z}$ for some $\beta > 0$. Then, there exists $u \in L^{\infty}([0,T];\mathbb{R}^m)$ such that
\begin{eqnarray}
&& c(t) + \sum_{i \in I_0(t)} u_i(t)a_i(t) = 0 ~ \mbox{a.e. in} ~ [0,T], \label{KKT_PL_2} \\
&& u_i(t) \geq 0  ~ \mbox{a.e. in} ~ [0,T], ~ i \in I. \label{NNeg_PL_2}
\end{eqnarray}
\end{theorem}

\begin{proof}
We keep the same notations as in the proof of Theorem \ref{CondNecPRL}. Assume that system below has a solution $(\alpha,\gamma) \in L^{\infty}([0,T];\mathbb{R} \times \mathbb{R}^n)$.
\begin{equation} \label{sistema22}
\left\{ \begin{array}{l}
\int_0^T c(t)^\top \gamma(t) dt < 0, \\ \delta_i^\beta(t) \bar{g}_i(t) \alpha(t) + \delta_i^\beta(t) a_i(t)^\top \gamma(t) \leq 0, ~ i \in I, ~ \mbox{a.e. in} ~ [0,T].
\end{array} \right.
\end{equation}
Let $K_1,K_2,K_3$ be positive scalars such that $\| \gamma(t) \| \leq K_1$, $\| a_i(t) \| \leq K_2$, $i \in I$, and $\| \alpha(t) \| \leq K_3$ a.e. in $[0,T]$. Let $0< \tau < \min\{1,\beta K_1^{-1}K_2^{-1},K_3^{-1}\}$. As in the proof of Theorem \ref{CondNecPRL}, we have that $F(\bar{z}+\tau\gamma) < F(\bar{z})$ and $a_i(t)^\top [\bar{z}(t) + \tau \gamma(t)] < b_i(t)$ a.e. in $[0,T]$ for $i \in I \setminus I_\beta(t)$ a.e. in $[0,T]$. For $i \in I_\beta(t)$, we obtain from the second inequality in (\ref{sistema22}) that $\bar{g}_i(t)\alpha(t) + a_i(t)^\top \gamma(t) \leq 0$ and
\begin{eqnarray*}
a_i(t)^\top [\bar{z}(t) + \tau \gamma(t)] - b_i(t) & = & [a_i(t)^\top \bar{z}(t) - b_i(t)] + \tau a_i(t)^\top \gamma(t) \\
& = & (1 - \tau\alpha(t))\bar{g}_i(t) \\ && + \tau [\bar{g}_i(t)\alpha(t) + a_i(t)^\top \gamma(t)] \leq 0.
\end{eqnarray*}
Thence, $a_i(t)^\top [\bar{z}(t) + \tau \gamma(t)] \leq b_i(t)$ a.e. in $[0,T]$. 

Summing up, we see that $\bar{z} + \tau \gamma \in \Omega$ with $F(\bar{z}+\tau\gamma) < F(\bar{z})$ for all $\tau > 0$ small enough, which is a contradiction to the optimality of $\bar{z}$ in (CLP). In conclusion, system (\ref{sistema22}) does not have any solution.

For almost every $t \in [0,T]$, we denote 
\begin{eqnarray*}
&& \bar{g}(t)=(\bar{g}_1(t),\ldots,\bar{g}_m(t)), \\
&& \Delta^\beta(t) = \mathrm{diag}(\delta_1^\beta(t),\ldots,\delta_m^\beta(t)), \\
&&M(t) = \left[ \begin{array}{ccc} \Delta^\beta(t) \bar{g}(t) & & \Delta^\beta(t) A(t) \end{array} \right] \in \mathbb{R}^{m \times (1+n)}.
\end{eqnarray*} 
Let us consider the continuous-time optimization problem below. 
\begin{equation} \label{problema_auxiliar12}
\begin{array}{ll} 
\mbox{Minimize} & \tilde{F}(\alpha,\gamma) = \displaystyle\int_{0}^{T} c(t)^\top \gamma(t) dt \\ 
\mbox{subject to} & M(t) \left[ \begin{array}{c} \alpha(t) \\ \gamma(t) \end{array} \right] \leq 0 ~ \mbox{a.e. in} ~ [0,T], \\
& (\alpha,\gamma) \in L^{\infty}([0,T];\mathbb{R} \times \mathbb{R}^n).
\end{array}
\end{equation}
It is easy to see that $(\bar{\alpha},\bar{\gamma})=(0,0)$ is an optimal solution. Indeed, it follows directly from that fact that system (\ref{sistema22}) does not have any solution.

We now consider the auxiliary unconstrained continuous-time optimization problem posed as follows:
\begin{equation} \label{problema_auxiliar22}
\begin{array}{ll} 
\mbox{minimize} & \Phi(\alpha,\gamma) = \displaystyle\int_{0}^{T} \left[ \begin{array}{c} 0 \\ c(t) \end{array} \right]^\top \left( \left[ \begin{array}{c} \alpha(t) \\ \gamma(t) \end{array} \right] - M(t)^+M(t)\left[ \begin{array}{c} \alpha(t) \\ \gamma(t) \end{array} \right] \right) dt \\ 
\mbox{subject to} & (\alpha,\gamma) \in L^{\infty}([0,T];\mathbb{R} \times \mathbb{R}^n).
\end{array}
\end{equation}
Set 
$$
\left[ \begin{array}{c} \hat{\alpha}(t) \\ \hat{\gamma}(t) \end{array} \right] = \left[ \begin{array}{c} \alpha(t) \\ \gamma(t) \end{array} \right] - M(t)^+M(t)\left[ \begin{array}{c} \alpha(t) \\ \gamma(t) \end{array} \right],
$$
where $(\alpha,\gamma) \in L^{\infty}([0,T];\mathbb{R} \times \mathbb{R}^n)$ is arbitrary. We have that $(\hat{\alpha},\hat{\gamma}) \in L^{\infty}([0,T];\mathbb{R} \times \mathbb{R}^n)$ and
$$
M(t) \left[ \begin{array}{c} \hat{\alpha}(t) \\ \hat{\gamma}(t) \end{array} \right] = 0 ~ \mbox{a.e. in} ~ [0,T],
$$
so that $(\hat{\alpha},\hat{\gamma})$ is a feasible solution to problem (\ref{problema_auxiliar12}). Then, $\tilde{F}(\hat{\alpha},\hat{\gamma}) \geq \tilde{F}(\bar{\alpha},\bar{\gamma})$, from where we obtain
$$
\Phi(\alpha,\gamma) = \tilde{F}(\hat{\alpha},\hat{\gamma}) \geq \tilde{F}(\bar{\alpha},\bar{\gamma}) = \Phi(\bar{\alpha},\bar{\gamma}).
$$
Thus, $(\bar{\alpha},\bar{\gamma})$ is an optimal solution of problem (\ref{problema_auxiliar22}). It follows from Proposition \ref{prob-sem-restr} that
$$
\left[ \begin{array}{c} 0 \\ c(t) \end{array} \right] - M(t)^\top (M(t)^+)^\top \left[ \begin{array}{c} 0 \\ c(t) \end{array} \right] = \left[ \begin{array}{c} 0 \\ 0 \end{array} \right] ~ \mbox{a.e. in} ~ [0,T].
$$
Defining
$$
u(t) = -(M(t)^+)^\top \left[ \begin{array}{c} 0 \\ c(t) \end{array} \right] ~ \mbox{a.e. in} ~ [0,T],
$$
we have that
\begin{equation} \label{eqt12}
\left[ \begin{array}{c} 0 \\ c(t) \end{array} \right] + \left[ \begin{array}{c} \bar{g}(t)^\top \Delta^\beta(t) \\ A(t)^\top \Delta^\beta(t) \end{array} \right] u(t) = \left[ \begin{array}{c} 0 \\ 0 \end{array} \right] ~ \mbox{a.e. in} ~ [0,T].
\end{equation}
Moreover, $u \in L^{\infty}([0,T];\mathbb{R}^m)$. Indeed, the norm of $M(t)^+$ is given by $\sigma_{\min}(M(t))^{-1}$, and this singular value is exactly the square root of the smallest positive eigenvalue of $M(t)M(t)^\top$. By noting that
$$
M(t)M(t)^\top = \Delta^\beta(t) \bar{g}(t) \bar{g}(t)^\top \Delta^\beta(t) + \Delta^\beta(t) A(t) A(t)^\top \Delta^\beta(t)
$$
and that $\mathrm{rank}(M(t)) = \mathrm{rank}(\Delta^\beta(t)A(t))$, bearing in mind the same reasoning as in the proof of Lemma \ref{Lema2}, it follows from the Weyl Theorem (see Theorem 4.3.1 in Horn and Johnson \cite{Horn:2013}, for instance) that
$$
\sigma_{\min}(M(t))^2 \geq \lambda_{1}(\Delta^\beta(t) \bar{g}(t) \bar{g}(t)^\top \Delta^\beta(t)) + \sigma_{\min}(\Delta^\beta(t) A(t))^2 \geq \sigma_{\min}(\Delta^\beta(t) A(t))^2.
$$
Hence, by using Lemma \ref{Lema1},
$$
\sigma_{\min}(M(t))^{-1} \leq \sigma_{\min}(\Delta^\beta(t) A(t))^{-1} = \sigma_{\min}(A^{I_\beta(t)}(t))^{-1} \leq C^{-1} ~ \mbox{a.e. in} ~ [0,T].
$$
Therefore, if $K$ is such that $|c(t)| \leq K$ a.e. in $[0,T]$,
$$
\| u(t) \| \leq \| (M(t)^+)^\top \| \| c(t) \| = \| M(t)^+ \| \|c(t) \| \leq C^{-1}K ~ \mbox{a.e. in} ~ [0,T].
$$

Let us show that $u_i(t) \geq 0$ for $i \in I_\beta(t)$ a.e. in $[0,T]$. Suppose that there exists $S \subset [0,T]$ with positive measure such that $u_{i_0}(t) < 0$ for some $i_0 \in I_\beta(t)$, for all $t \in S$. Define, almost everywhere in $[0,T]$, $d(t) \in \mathbb{R}^{m}$ as $d_i(t) = 0$ for $i \in I$, $i \neq i_0$,
$$
d_{i_0}(t) = \left\{ \begin{array}{rl}
-1, & \mbox{if} ~ t \in S, \\ 0, & \mbox{otherwise}.
\end{array} \right.
$$
Clearly, $d \in L^{\infty}([0,T];\mathbb{R}^m)$. We know from the hypothesis that 
matrix $A^{I_\beta(t)}(t)A^{I_\beta(t)}(t)^\top$ is invertible a.e. in $[0,T]$. Take $\tilde{\alpha}(t) = 1$ and
$$
\tilde{\gamma}(t) = A^{I_\beta(t)}(t)^\top [A^{I_\beta(t)}(t) A^{I_\beta(t)}(t)^\top]^{-1} d^{I_\beta(t)}(t) ~ \mbox{a.e. in} ~ [0,T].
$$
It is clear that $(\tilde{\alpha},\tilde{\gamma}) \in L^{\infty}([0,T];\mathbb{R} \times \mathbb{R}^n)$. We have also that
$$
A^{I_\beta(t)}(t) \tilde{\gamma}(t) = d^{I_\beta(t)}(t) \leq 0 ~ \mbox{a.e. in} ~ [0,T].
$$
Therefore,
$$
M(t) \left[ \begin{array}{c} \tilde{\alpha}(t) \\ \tilde{\gamma}(t) \end{array} \right] \leq 0 ~ \mbox{a.e. in} ~ [0,T],
$$
that is, $(\tilde{\alpha},\tilde{\gamma})$ is a feasible solution to (\ref{problema_auxiliar12}). It follows that $\tilde{F}(\tilde{\alpha},\tilde{\gamma}) \geq \tilde{F}(\bar{\alpha},\bar{\gamma}) = 0$. On the other hand, by using (\ref{eqt12}),
\begin{eqnarray*}
\tilde{F}(\tilde{\alpha},\tilde{\gamma}) &=& \int_0^T c(t)^\top \tilde{\gamma}(t) dt = -\int_0^T u(t)^\top \Delta^\beta(t) A(t) \tilde{\gamma}(t) dt \\ & = &  -\int_0^T u(t)^\top d(t) dt = -\int_0^T u_{i_0}(t)d_{i_0}(t) dt = \int_S u_{i_0}(t) dt < 0,
\end{eqnarray*}
a contradiction.

From (\ref{eqt12}), it comes
\begin{eqnarray*}
&& \sum_{i \in I_\beta(t)} u_i(t) \bar{g}_i(t) = 0  ~ \mbox{a.e. in} ~ [0,T], \\ 
&& c(t) + \sum_{i \in I_\beta(t)} u_i(t)a_i(t) = 0 ~ \mbox{a.e. in} ~ [0,T].
\end{eqnarray*}
The first equality above together with $\bar{g}_i(t) \leq 0$ and $u_i(t) \geq 0$ for $i \in I_\beta(t)$ a.e. in $[0,T]$ imply in $u_i(t) = 0$ for $i \notin I_0(t)$ a.e. in $[0,T]$. Therefore, (\ref{KKT_PL_2})-(\ref{NNeg_PL_2}) are valid. 
\end{proof}

Given the linearity of the problem, the converse of Theorem \ref{CondNecPRL2} holds true.

\begin{theorem}
Let $\bar{z} \in \Omega$ and assume that there exists $u \in L^{\infty}([0,T];\mathbb{R}^m)$ such that
\begin{eqnarray}
&& c(t) + \sum_{i \in I_0(t)} u_i(t)a_i(t) = 0 ~ \mbox{a.e. in} ~ [0,T], \label{KKT_PL_3} \\
&& u_i(t) \geq 0  ~ \mbox{a.e. in} ~ [0,T], ~ i \in I. \label{NNeg_PL_3}
\end{eqnarray} 
Then $\bar{z}$ is an optimal solution of (CLP).
\end{theorem}

\begin{proof}
Let $z \in \Omega$ be an arbitrary feasible solution of (CLP). Then $A(t)z(t) \leq b(t)$ a.e. in $[0,T]$, that is, $a_i(t)^\top z(t) \leq b_i(t)$ a.e. in $[0,T]$, $i \in I$. For $i \in I_0(t)$, we have $a_i(t)^\top \bar{z}(t) = b_i(t)$ a.e. in $[0,T]$. Hence,
\begin{equation} \label{Desig1}
a_i(t)^\top [z(t)-\bar{z}(t)] \leq 0, ~ i \in I_0(t) ~ \mbox{a.e. in} ~ [0,T].
\end{equation}
By (\ref{KKT_PL_3}), (\ref{NNeg_PL_3}) and (\ref{Desig1}), we obtain
\begin{eqnarray*}
F(z) - F(\bar{z}) & = & \int_0^T c(t)^\top [z(t)-\bar{z}(t)] dt \\
& = & -\int_0^T \sum_{i \in I_0(t)} u_i(t)a_i(t)^\top [z(t)-\bar{z}(t)] dt \geq 0.
\end{eqnarray*}
It follows that $\bar{z}$ is an optimal solution of (CLP). 
\end{proof}

In what follows, we present a simple academic example to illustrate the results above.

\begin{example} \label{Exemplo}
Let us consider the following linear continuous-time optimization problem:
$$
\begin{array}{ll} 
\mbox{minimize} & F(z) = \displaystyle\int_{0}^{2} c(t)^\top z(t) dt \\ 
\mbox{subject to} & A(t)z(t) \leq b(t) ~ \mbox{a.e. in} ~ [0,2], \\
& z \in L^{\infty}([0,2];\mathbb{R}^2),
\end{array}
$$
where, for almost every $t \in [0,1]$,
$$
c(t) = \left[ \begin{array}{c} t-1 \\ -1 \end{array} \right], \quad
A(t) = \left[ \begin{array}{rr} 0 & -1 \\ -1 & 0 \\ -1 & 1\\ 1 & 1 \\ 0 & 1 \end{array} \right], \quad
b(t) = \left[ \begin{array}{c} 0 \\ 0 \\ 0 \\ 3 \\ \frac{1}{4}+\frac{5}{8}t \end{array} \right]
$$
and, for almost every $t \in (1,2]$,
$$
c(t) = \left[ \begin{array}{c} -t+1 \\ -1 \end{array} \right], \quad
A(t) = \left[ \begin{array}{rr} 0 & -1 \\ -1 & 0 \\ 1 & -1 \\ 1 & 1 \\ 0 & 1 \end{array} \right], \quad
b(t) = \left[ \begin{array}{c} 0 \\ 0 \\ 0 \\ 3 \\ \frac{1}{4}+\frac{5}{8}t \end{array} \right].
$$
Below we have illustrations of the feasible region for the values of $t=0$, $t=1/2$, $t=1$, $t=3/2$ and $t=2$:
\begin{center}
\includegraphics[scale=0.28]{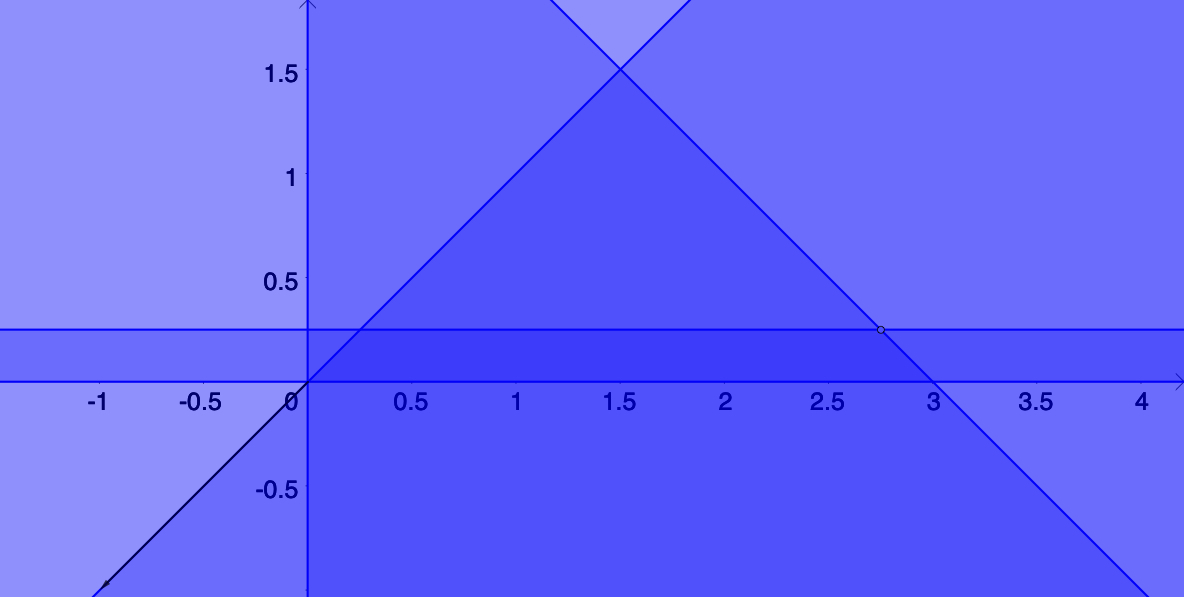}
\end{center}
\begin{center}
\includegraphics[scale=0.28]{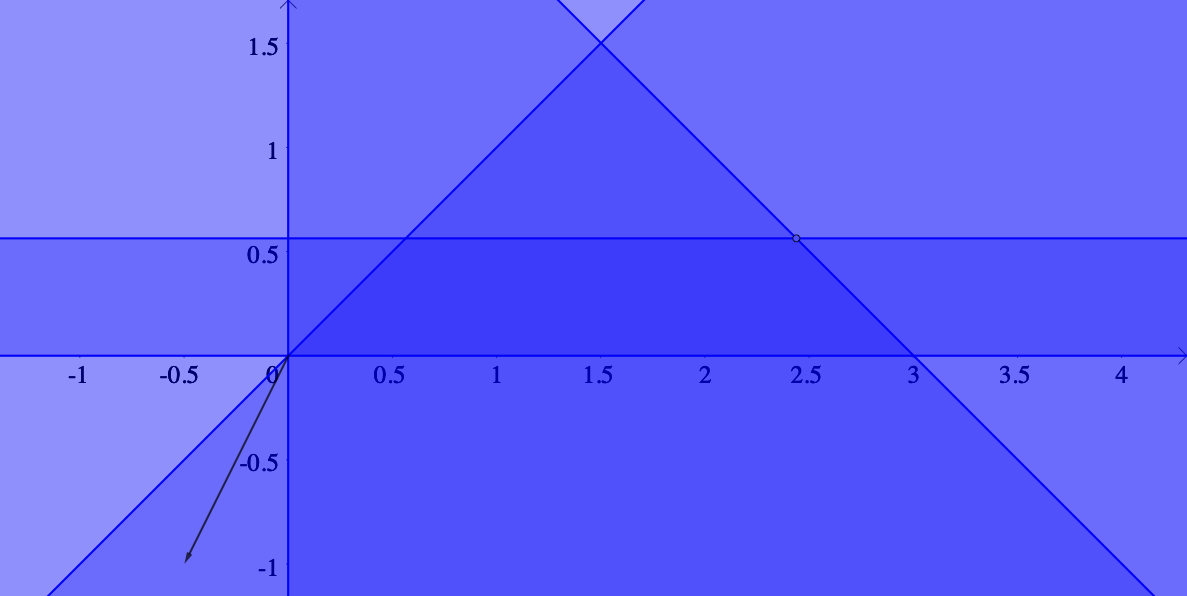}
\end{center}
\begin{center}
\includegraphics[scale=0.28]{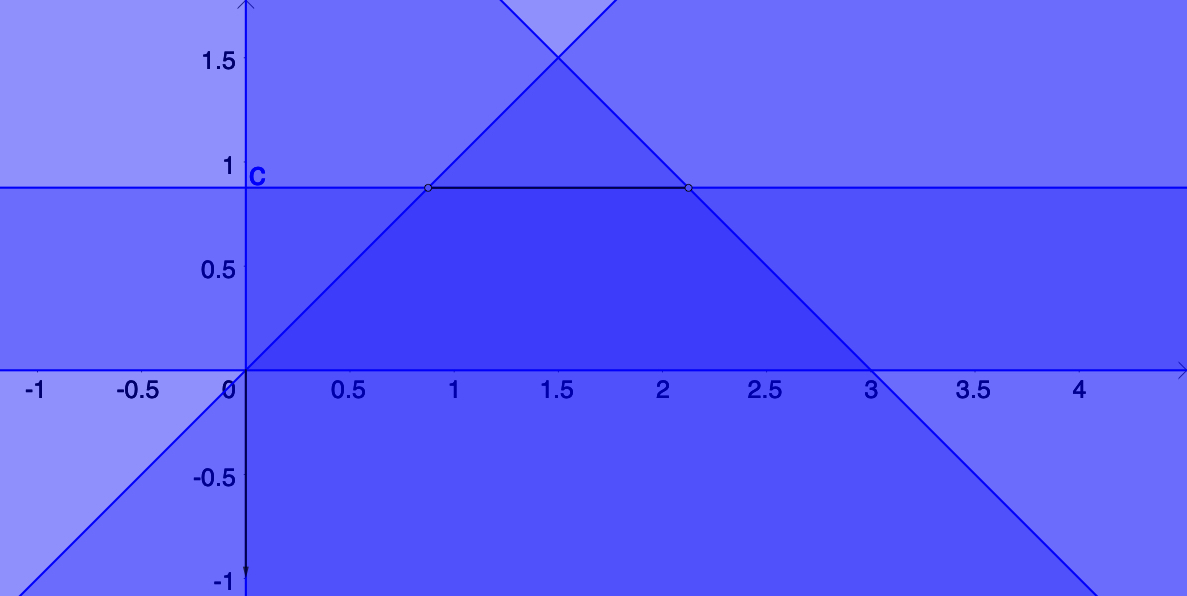}
\end{center}
\begin{center}
\includegraphics[scale=0.28]{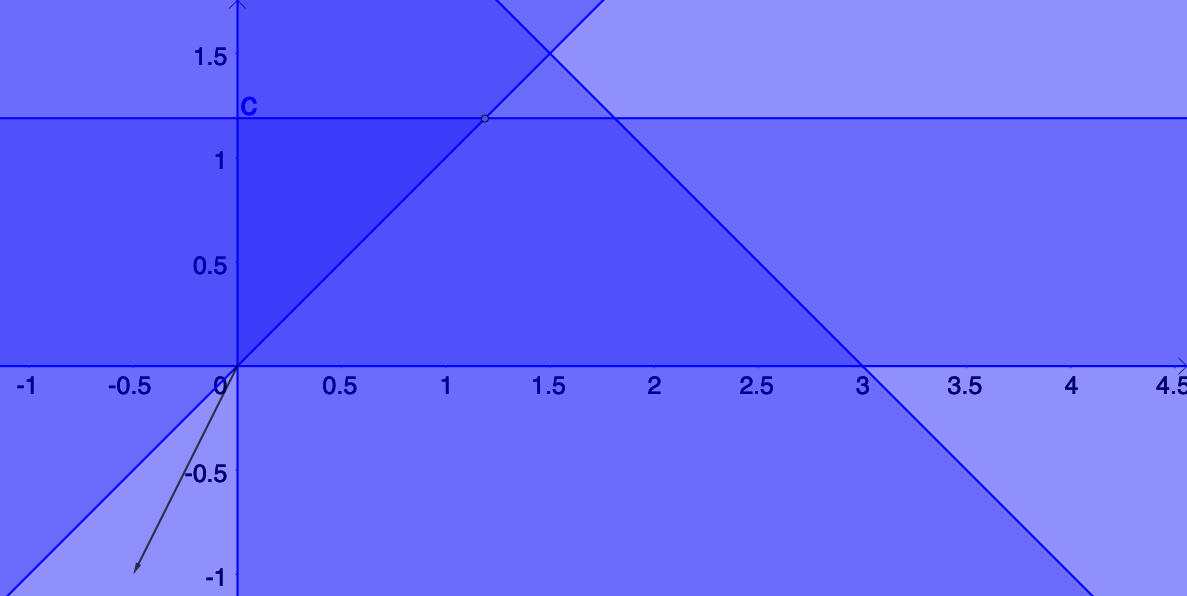}
\end{center}
\begin{center}
\includegraphics[scale=0.28]{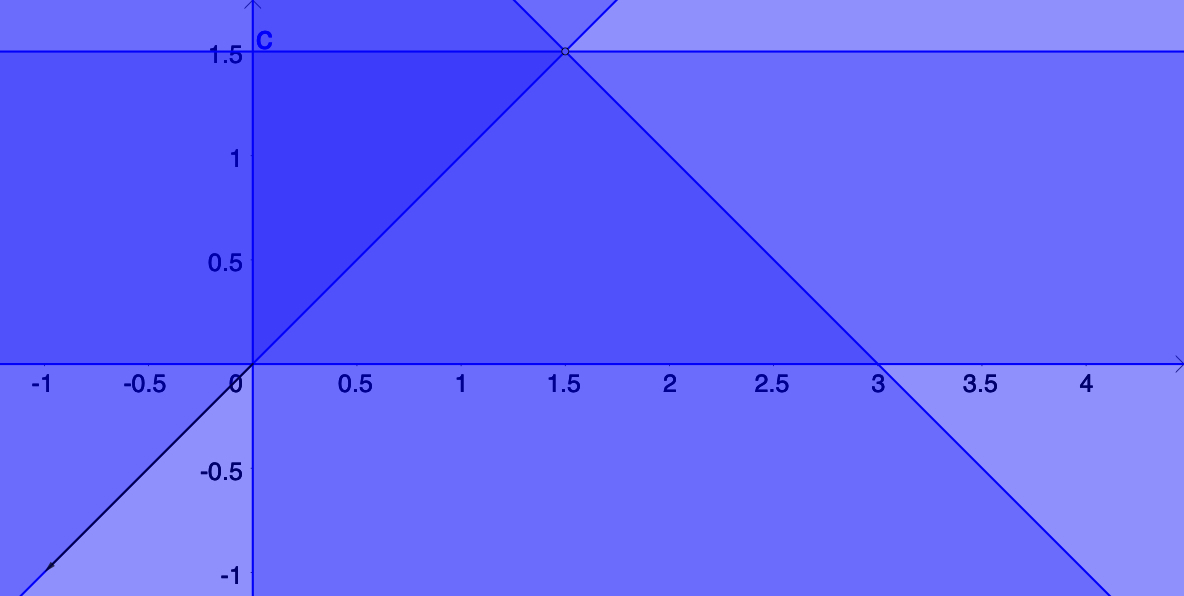}
\end{center}
Let $\bar{z} : [0,2] \rightarrow \mathbb{R}^2$ be given as
$$
\bar{z}_1(t) = \left\{ \begin{array}{l}
\frac{11}{4}-\frac{5}{8}t ~ \mbox{a.e. in} ~ [0,1], \\
\frac{1}{4}+\frac{5}{8}t ~ \mbox{a.e. in} ~ (1,2],
\end{array} \right.
$$
and
$$
\bar{z}_2(t) = \left\{ \begin{array}{ll}
\frac{1}{4}+\frac{5}{8}t ~ \mbox{a.e. in} ~ [0,1], \\
\frac{1}{4}+\frac{5}{8}t ~ \mbox{a.e. in} ~ (1,2].
\end{array} \right.
$$
For almost every $t \in [0,1]$, we have
$$
c(t)^\top \bar{z}(t) = \left[ \begin{array}{cc} t-1 & -1 \end{array} \right] \left[ \begin{array}{c} \frac{11}{4}-\frac{5}{8}t \\ \frac{1}{4}+\frac{5}{8}t \end{array} \right] = -\frac{5}{8}t^2+\frac{11}{4}t-3.
$$
For almost every $t \in (1,2]$,
$$
c(t)^\top \bar{z}(t) = \left[ \begin{array}{cc} -t+1 & -1 \end{array} \right] \left[ \begin{array}{c} \frac{1}{4}+\frac{5}{8}t \\ \frac{1}{4}+\frac{5}{8}t \end{array} \right] = -\frac{5}{8}t^2-\frac{1}{4}t.
$$
Then,
\begin{eqnarray*}
F(\bar{z}) &=& \int_0^1 \left( -\frac{5}{8}t^2+\frac{11}{4}t-3 \right) dt + \int_1^2 \left( -\frac{5}{8}t^2-\frac{1}{4}t \right) dt \\
&=& -\frac{11}{6} - \frac{11}{6} = -\frac{11}{3}.  
\end{eqnarray*}
Take $\beta = \frac{1}{8}$. Then, for almost every $t \in [0,1]$,
\begin{eqnarray*}
&& a_1(t)^\top \bar{z}(t) - b_1(t) = -\frac{1}{4}-\frac{5}{8}t \leq -\frac{1}{4} < -\beta, \\
&& a_2(t)^\top \bar{z}(t) - b_2(t) = -\frac{11}{4}+\frac{5}{8}t \leq -\frac{17}{8} < -\beta, \\
&& a_3(t)^\top \bar{z}(t) - b_3(t) = -\frac{10}{4}-\frac{5}{4}t \leq -\frac{10}{4} < -\beta, \\
&& a_4(t)^\top \bar{z}(t) - b_4(t) = 0 \geq -\beta, \\
&& a_5(t)^\top \bar{z}(t) - b_5(t) = 0 \geq -\beta,
\end{eqnarray*}
and, for almost every $t \in (1,2]$,
\begin{eqnarray*}
&& a_1(t)^\top \bar{z}(t) - b_1(t) = -\frac{1}{4}-\frac{5}{8}t \leq -\frac{1}{4} < -\beta, \\
&& a_2(t)^\top \bar{z}(t) - b_2(t) = -\frac{1}{4}-\frac{5}{8}t \leq -\frac{1}{4} < -\beta, \\
&& a_3(t)^\top \bar{z}(t) - b_3(t) = 0 \geq -\beta, \\
&& a_4(t)^\top \bar{z}(t) - b_4(t) = -\frac{5}{2}+\frac{5}{4}t \\
&& a_5(t)^\top \bar{z}(t) - b_5(t) = 0 \geq -\beta.
\end{eqnarray*}
Let us analyze the constraint $a_4(t)^\top \bar{z}(t) - b_4(t) \leq 0$ in this second case. We have that $a_4(t)^\top \bar{z}(t) - b_4(t) < -\beta$ a.e. in $t \in (1,19/10)$ and $a_4(t)^\top \bar{z}(t) - b_4(t) \geq -\beta$ a.e. in $t \in (19/10,2]$. Thus,
$$
I_0(t) = \left\{ \begin{array}{l}
\{4,5\} ~ \mbox{a.e. in} ~ [0,1], \\
\{3,5\} ~ \mbox{a.e. in} ~ (1,2), \\
\{3,4,5\} ~ \mbox{se} ~ t=2,
\end{array} \right.
$$
and
$$
I_\beta(t) = \left\{ \begin{array}{l}
\{4,5\} ~ \mbox{a.e. in} ~ [0,1], \\
\{3,5\} ~ \mbox{a.e. in} ~ (1,\frac{19}{10}), \\
\{3,4,5\} ~ \mbox{a.e. in} ~ [\frac{19}{10},2].
\end{array} \right.
$$
Taking $I_\sharp(t) = \{4,5\}$ a.e. in $[0,1]$ and $I_\sharp(t) = \{3,5\}$ a.e. in $(1,2]$, it follows that $\beta$-CR is satisfied at $\bar{z}$ (note that (i) $\beta$-FR is not satisfied at $\bar{z}$; (ii) $I_\sharp(t) \subset I_0(t)$ a.e. in $[0,T]$). It is clear that (A1)-(A2) are valid. We will show that $\bar{z}$ is an optimal solution. The necessary optimality condition is satisfied with $u_1(t) = u_2(t) = 0$ a.e. in $[0,2]$,
$$
u_3(t) = \left\{ \begin{array}{l} 0 ~ \mbox{a.e. in} ~ [0,1], \\ t-1 ~ \mbox{a.e. in} ~ (1,2], \end{array} \right. \quad
u_4(t) = \left\{ \begin{array}{l} -t+1 ~ \mbox{a.e. in} ~ [0,1], \\ 0 ~ \mbox{a.e. in} ~ (1,2], \end{array} \right.
$$
and $u_5(t) = t$ a.e. in $[0,2]$. Provided the problem is linear, we conclude that $\bar{z}$ is an optimal solution. Below we have an illustration of the optimal solution $\bar{z}$:
\begin{center}
\includegraphics[scale=0.28]{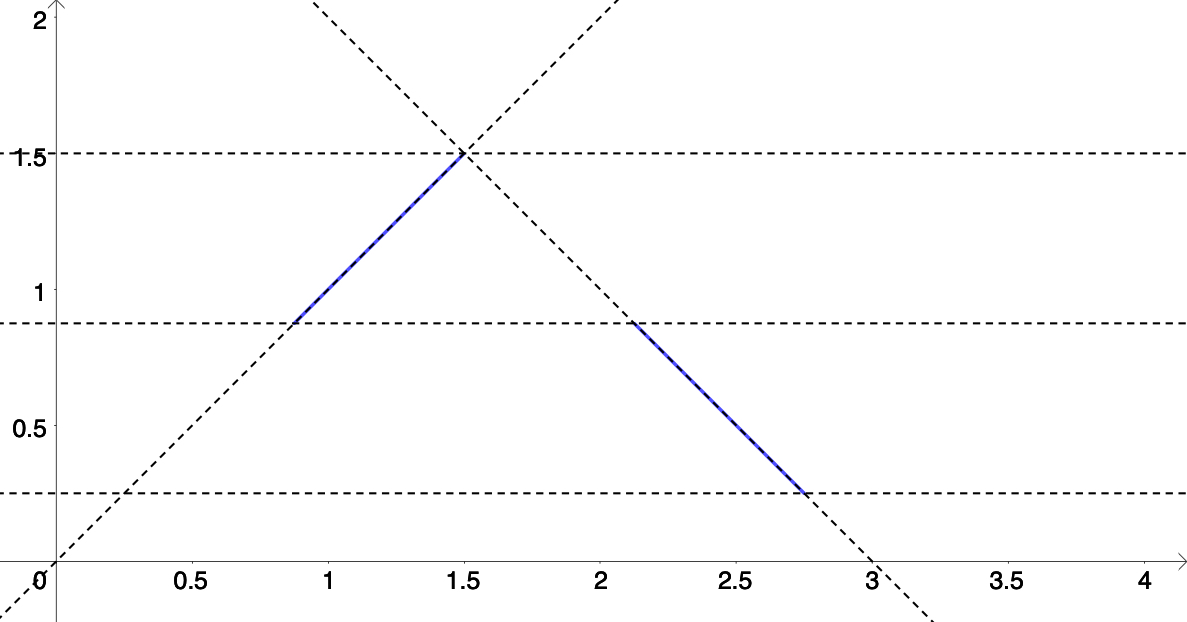}
\end{center} 
\end{example}

\section{Duality results} \label{dualidade}

Here we give the definition of the dual problem and provide the primal-dual relationship. Among others, weak and strong duality properties and the complementary slackness theorem are presented.

Associated to (CLP), we define the following dual problem:
$$
\begin{array}{ll} 
\mbox{maximize} & G(w) = \displaystyle\int_{0}^{T} b(t)^\top w(t) dt \\ 
\mbox{subject to} & A(t)^\top w(t) = c(t) ~ \mbox{a.e. in} ~ [0,T], \\
& w(t) \leq 0 ~ \mbox{a.e. in} ~ [0,T], \\
& w \in L^{\infty}([0,T];\mathbb{R}^m).
\end{array}
\eqno{\mathrm{(CDP)}}
$$

From now on, (CLP) will be referred to as the primal problem. 

The set of all feasible dual solutions will be denoted by $\Theta$, i.e.,
$$
\Theta = \{ w \in L^{\infty}([0,T];\mathbb{R}^m) : A(t)^\top w(t) = c(t), ~ w(t) \leq 0 ~ \mbox{a.e. in} ~ [0,T] \}.
$$
We say that $\bar{w} \in \Theta$ is an optimal solution of (CDP) if
$$
G(\bar{w}) \geq G(w) ~ \forall w \in \Theta.
$$

We begin with the weak duality property. 

\begin{theorem}[Weak Duality] \label{DualidadeFraca}
Let $z \in \Omega$ and $w \in \Theta$ be feasible solutions of (CLP) and (CDP). Then, $G(w) \leq F(z)$.
\end{theorem}

We will omit the proof of the last result since it is standard and straightforward. The same applies to the two next ones.


\begin{theorem}
Assume that $\Omega \neq \emptyset$. If there exists a sequence $\{ z^k \} \subset \Omega$ such that $F(z^k) \to -\infty$ when $k \to \infty$, then $\Theta = \emptyset$.
\end{theorem}


\begin{theorem}
Assume that $\Theta \neq \emptyset$. If there exists a sequence $\{ w^k \} \subset \Theta$ such that $G(w^k) \to +\infty$ when $k \to \infty$, then $\Omega = \emptyset$.
\end{theorem}


\begin{lemma} \label{LemaDual}
Let $\bar{z} \in \Omega$ be an optimal solution of (CLP). Assume that $\beta$-RC is satisfied at $\bar{z}$ for some $\beta > 0$. Then $\tilde{w} = -u$, where $u \in L^{\infty}([0,T];\mathbb{R}^n)$ is the Lagrange multiplier associated to $\bar{z}$, given in Theorem \ref{CondNecPRL}, is a feasible solution of (CDP), that is, $\tilde{w} \in \Theta$.
\end{lemma}

\begin{proof}
From (\ref{KKT_PL}), we have $c(t) = -A(t)^\top u(t)$ a.e. in $[0,T]$, so that $A(t)^\top \tilde{w}(t) = c(t)$ a.e. in $[0,T]$. From (\ref{NNeg_PL}), it follows that $\tilde{w}(t) = -u(t) \leq 0$ a.e. in $[0,T]$. Therefore, $\tilde{w} \in \Theta$. 
\end{proof}

Analogously, we have

\begin{lemma} \label{LemaDual2}
Let $\bar{z} \in \Omega$ be an optimal solution of (CLP). Assume that $\beta$-FR is satisfied at $\bar{z}$ for some $\beta > 0$. Then $\tilde{w} = -u$, where $u \in L^{\infty}([0,T];\mathbb{R}^n)$ is the Lagrange multiplier associated to $\bar{z}$, given in Theorem \ref{CondNecPRL2}, is a feasible solution of (CDP), that is, $\tilde{w} \in \Theta$.
\end{lemma}

\begin{theorem} \label{RelPrimalDual1}
Assume that  either $\beta$-RC or $\beta$-FR is satisfied at each $z \in \Omega$ for some $\beta=\beta(z) > 0$. In addition, assume that $\Omega \neq \emptyset$. If $\Theta = \emptyset$, then (CLP) does not have any optimal solution.
\end{theorem}

\begin{proof}
If (CLP) has an optimal solution, by the last two lemmas, there exists $\tilde{w} \in \Theta$, contradicting the hypothesis. 
\end{proof}

The hypothesis below will be necessary in what follows.

\begin{description}
\item[(H)] The hypothesis H is said to be satisfied if for each $\bar{w} \in \Theta$ there exists a constant $K_{\bar{w}} > 0$ such that $\det(\Upsilon(t)\Upsilon(t)^\top) \geq K_{\bar{w}}$ a.e. in $[0,T]$ where
$$
\Upsilon(t) = \left[
\begin{array}{cc}
A(t)^\top & 0 \\ -I & D(t)
\end{array}
\right],
$$
with $D(t) := \mathrm{diag}(\{ -2 \sqrt{-\bar{w}_i(t)} \}_{i \in I} )$ a.e. in $[0,T]$ and $I$ denotes the identity matrix of order $m$.
\end{description}

\begin{theorem} \label{Teorema9}
Assume that hypothesis H is satisfied. In addition, assume that $\Theta \neq \emptyset$. If $\Omega = \emptyset$, then (CDP) does not have any optimal solution.
\end{theorem}

\begin{proof}
Suppose that (CDP) has an optimal solution, say $\bar{w} \in \Theta$. Let us note that hypothesis H implies that the full rank assumption in Monte and de Oliveira \cite{monte:2019} is satisfied. It follows from Theorem 4.2 in \cite{monte:2019}, that there exist $u \in L^\infty([0,T];\mathbb{R}^n)$ and $v \in L^\infty([0,T];\mathbb{R}^m)$ such that $b(t) + A(t)u(t) - v(t) = 0$, $v(t) \geq 0$ and $v(t)^\top \bar{w}(t) = 0$ a.e. in $[0,T]$. Taking $\tilde{z} = -u$, we have $A(t)\tilde{z}(t) = b(t) - v(t) \leq b(t)$ a.e. in $[0,T]$, that is, $\tilde{z} \in \Omega$, contradicting the hypothesis.
\end{proof}

In the following results, we will make use of the assumption stated next. Let $\beta > 0$. 

\begin{description}
\item[($\beta$-A)] The assumption $\beta$-A is said to be satisfied at $\bar{z} \in \Omega$ if either
\begin{enumerate}
\item[(i)] $\beta$-RC is valid at $\bar{z}$ with $I_\sharp(t) \subseteq I_0(t)$ a.e. in $[0,T]$; or
\item[(ii)] $\beta$-FR is valid at $\bar{z}$.
\end{enumerate}
\end{description}

\begin{theorem}[Strong Duality] \label{DualidadeForte}
Let $\bar{z} \in \Omega$ and $\bar{w} \in \Theta$ be feasible solutions of (CLP) and (CDP). Assume that $\beta$-A is satisfied at $\bar{z}$ for some $\beta > 0$. Then $F(\bar{z}) = G(\bar{w})$ if, and only if, $\bar{z}$ is an optimal solution of (CLP) and $\bar{w}$ is an optimal solution of (CDP).
\end{theorem}

\begin{proof}
Suppose that $F(\bar{z}) = G(\bar{w})$. It follows from the weak duality property, Theorem \ref{DualidadeFraca}, that $F(z) \geq G(\bar{w}) = F(\bar{z})$ for all $z \in \Omega$. Thus, $\bar{z}$ is an optimal solution of (CLP). It follows from the same property that $G(w) \leq F(\bar{z}) = G(\bar{w})$ for all $w \in \Theta$, so that $\bar{w}$ is an optimal solution of (CDP).

Reciprocally, suppose that $\bar{z}$ is an optimal solution of (CLP) and $\bar{w}$ is an optimal solution of (CDP). Since $\beta$-A is valid at $\bar{z}$, either (i) $\beta$-CR is valid with $I_\sharp(t) \subseteq I_0(t)$ a.e. in [0,T] or (ii) $\beta$-FR is satisfied. We assume first that (ii) holds.  We know from Lemma \ref{LemaDual2} that $\tilde{w} = -u \in \Theta$, where $u$ is the Lagrange multiplier associated to $\bar{z}$, given in Theorem \ref{CondNecPRL2}. We know from the proof of Theorem \ref{CondNecPRL2} that $u_i(t) = 0$ for $i \notin I_0(t)$ a.e. in $[0,T]$. Therefore, from (\ref{KKT_PL_2}) we have that $A^{I_0(t)}(t)^\top u^{I_0(t)}(t) = A(t)^\top u(t) = -c(t)$ a.e. in $[0,T]$. For $i \in I_0(t)$ we have $a_i(t)^\top \bar{z}(t) = b_i(t)$ a.e. in $[0,T]$. In vectorial notation, $A^{I_0(t)}(t)\bar{z}(t) = b^{I_0(t)}$ a.e. in $[0,T]$. Keeping this in mind, we obtain
\begin{eqnarray*}
G(\bar{w}) & \geq & G(\tilde{w}) = \int_{0}^{T} b(t)^\top \tilde{w}(t) dt = - \int_{0}^{T} b(t)^\top u(t) dt \\
& = & - \int_{0}^{T} b^{I_0(t)}(t)^\top u^{I_0(t)}(t) dt = - \int_{0}^{T} \bar{z}(t)^\top A^{I_0(t)}(t)^\top u^{I_0(t)}(t) dt \\
& = & \int_0^T \bar{z}(t)^\top c(t) dt = F(\bar{z}).
\end{eqnarray*}
On the other hand, we have from the weak duality property that $G(\bar{w}) \leq F(\bar{z})$. Thus, $G(\bar{w}) = F(\bar{z})$. When (ii) holds, we use Lemma \ref{LemaDual} and Theorem \ref{CondNecPRL} instead of Lemma \ref{LemaDual2} and Theorem \ref{CondNecPRL2}, and the proof follows analogously by noting that $I_\sharp(t) \subseteq I_0(t)$ implies in $u_i(t) = 0$ for $i \notin I_0(t)$ a.e. in $[0,T]$.
\end{proof}

Observe in the proof of the last theorem that the assumption that $\beta$-A is satisfied was only needed to prove that the optimal values of the primal and dual problems coincide. The other statement, which is a direct consequence of weak duality, always holds.

\begin{corollary} \label{CorolarioMultLag}
Let $\bar{z} \in \Omega$ be an optimal solution of (CLP). Assume that $\beta$-RC is satisfied at $\bar{z}$ with $I_\sharp(t) \subseteq I_0(t)$ a.e. in $[0,T]$ for some $\beta > 0$. Then $\tilde{w} = -u$, where $u \in L^{\infty}([0,T];\mathbb{R}^n)$ is the Lagrange multiplier associated to $\bar{z}$, given in Theorem \ref{CondNecPRL}, is an optimal solution of (CDP) with $G(\tilde{w}) = F(\bar{z})$.
\end{corollary}

\begin{proof}
We know from Lemma \ref{LemaDual} that $\tilde{w} \in \Theta$ and we saw in the proof of Theorem \ref{DualidadeForte} that $G(\tilde{w}) = F(\bar{z})$. It follows from Theorem \ref{DualidadeForte} that $\tilde{w}$ is an optimal solution of (CDP).
\end{proof}

Similarly,

\begin{corollary} \label{CorolarioMultLag2}
Let $\bar{z} \in \Omega$ be an optimal solution of (CLP). Assume that $\beta$-FR is satisfied at $\bar{z}$ for some $\beta > 0$. Then $\tilde{w} = -u$, where $u \in L^{\infty}([0,T];\mathbb{R}^n)$ is the Lagrange multiplier associated to $\bar{z}$, given in Theorem \ref{CondNecPRL2}, is an optimal solution of (CDP) with $G(\tilde{w}) = F(\bar{z})$.
\end{corollary}

\begin{example} \label{Exemplo2}
We now revisit Example \ref{Exemplo}. The dual problem is given as
$$
\begin{array}{ll} 
\mbox{maximize} & G(w) = \displaystyle\int_{0}^{2} \left[3w_4(t) + \left(\frac{1}{4}+\frac{5}{8}t\right)w_5(t)\right] dt \\ 
\mbox{subject to} & -w_2(t)-w_3(t)+w_4(t) = t-1 ~ \mbox{a.e. in} ~ [0,1], \\ 
& -w_2(t)+w_3(t)+w_4(t) = -t+1 ~ \mbox{a.e. in} ~ (1,2], \\
& -w_1(t)+w_3(t)+w_4(t)+w_5(t) = -1 ~ \mbox{a.e. in} ~ [0,1], \\
& -w_1(t)-w_3(t)+w_4(t)+w_5(t) = -1 ~ \mbox{a.e. in} ~ (1,2], \\
& w_i(t) \leq 0 ~ \mbox{a.e. in} ~ [0,T], ~ i = 1,2,3,4,5, \\
& w \in L^{\infty}([0,T];\mathbb{R}^m).
\end{array}
$$
We know from Example \ref{Exemplo} that $\beta$-RC is satisfied at the optimal solution $\bar{z}$ with $I_\sharp(t) \subseteq I_0(t)$ a.e. in $[0,2]$. It follows from Corollary \ref{CorolarioMultLag} that
\begin{align*}
& \tilde{w}_1(t) = \tilde{w}_2(t) = 0 ~ \mbox{a.e. in} ~  [0,2], \\
& \tilde{w}_3(t) = \left\{ \begin{array}{l} 0 ~ \mbox{a.e. in} ~ [0,1], \\ -t+1 ~ \mbox{a.e. in} ~ (1,2], \end{array} \right. \\
& \tilde{w}_4(t) = \left\{ \begin{array}{l} t-1 ~ \mbox{a.e. in} ~ [0,1], \\ 0 ~ \mbox{a.e. in} ~ (1,2], \end{array} \right. \\
& \tilde{w}_5(t) = -t ~ \mbox{a.e. in} [0,2],
\end{align*}
is an optimal solution of the dual problem above. We also have
\begin{eqnarray*}
G(\tilde{w}) &=& \int_0^1 \left[ 3(t-1) + \left( \frac{1}{4}+\frac{5}{8}t \right) (-t) \right] dt + \int_1^2 \left( \frac{1}{4}+\frac{5}{8}t \right) (-t) dt \\
&=& -\frac{11}{6} -\frac{11}{6} = -\frac{11}{3} = F(\bar{z}).
\end{eqnarray*}
\end{example}

The continuous-time linear problem in the example above, although it is simple, cannot be treated by some previous results. For example, the assumptions in the strong duality theorem given in Wu \cite{Wu:2013} are not satisfied: an assumption on the matrix defining the system of linear inequality constraints is not satisfied. Indeed, to see this, let us put (CLP) in the format considered in reference \cite{Wu:2013}. Let $z = z^+ - z^-$. (CLP) can be written as
$$
\begin{array}{ll} 
\mbox{maximize} & -F(z) = \displaystyle\int_{0}^{T} \left[ \begin{matrix} -c(t)^\top & c(t)^\top \end{matrix} \right] \left[ \begin{matrix} z^+(t) \\ z^-(t) \end{matrix} \right] dt \\ 
\mbox{subject to} & \left[ \begin{matrix} A(t) & -A(t) \end{matrix} \right] \left[ \begin{matrix} z^+(t) \\ z^-(t) \end{matrix} \right] \leq b(t) ~ \mbox{a.e. in} ~ [0,T], \\
& \left[ \begin{matrix} z^+(t) \\ z^-(t) \end{matrix} \right] \geq 0 ~ \mbox{a.e. in} ~ [0,T], \\
& z^+,z^- \in L^{\infty}([0,T];\mathbb{R}^n).
\end{array}
$$
In \cite{Wu:2013}, it is assumed that $\sum_{i=1}^{m}b_{ij}(t) > 0$ a.e. in $[0,T]$ for each $j$, where $B(t) = \left[ \begin{matrix} A(t) & -A(t) \end{matrix} \right]$ a.e. in $[0,T]$. This is not satisfied for the matrix $B(t)$ in example above (see the definition of the matrix $A(t)$ in Example \ref{Exemplo}).

\begin{theorem}[Complementary Slackness Theorem] \label{FolgasCompl}
Let $\bar{z} \in L^{\infty}([0,T];\mathbb{R}^n)$ and $\bar{w} \in L^{\infty}([0,T];\mathbb{R}^m)$. Assume that $\beta$-A is satisfied at $\bar{z}$ for some $\beta > 0$. Then, solutions $\bar{z} \in L^{\infty}([0,T];\mathbb{R}^n)$ and $\bar{w} \in L^{\infty}([0,T];\mathbb{R}^m)$ are optimal for (CLP) and (CDP) if, and only if, there exists $\bar{y} \in L^{\infty}([0,T];\mathbb{R}^m)$ such that
\begin{eqnarray}
&& A(t)\bar{z}(t) + \bar{y}(t) = b(t), ~ \bar{y}(t) \geq 0 ~ \mbox{a.e. in} ~ [0,T], \label{FC1} \\
&& A(t)^\top \bar{w}(t) = c(t), ~ \bar{w}(t) \leq 0 ~ \mbox{a.e. in} ~ [0,T], \label{FC2} \\
&& \bar{y}(t)^\top \bar{w}(t) = 0 ~ \mbox{a.e. in} ~ [0,T]. \label{FC3}
\end{eqnarray}
\end{theorem}

\begin{proof}
Suppose that $\bar{z} \in L^{\infty}([0,T];\mathbb{R}^n)$ and $\bar{w} \in L^{\infty}([0,T];\mathbb{R}^m)$ are optimal for (CLP) and (CDP). It is clear that $\bar{z} \in \Omega$ and that $\bar{w} \in \Theta$. Defining $\bar{y}(t) = b(t) - A(t)\bar{z}(t)$ a.e. in $[0,T]$ and taking into account that $\bar{z} \in \Omega$, it follows that (\ref{FC1}) holds. The validity of (\ref{FC2}) follows directly from $\bar{w} \in \Theta$. If $\tilde{w} = -u$, where $u \in L^{\infty}([0,T];\mathbb{R}^n)$ is the Lagrange multiplier associated to $\bar{z}$, we know from Corollaries \ref{CorolarioMultLag} and \ref{CorolarioMultLag2} that $\tilde{w}$ is an optimal solution of (CDP) with $G(\tilde{w}) = F(\bar{z})$. Therefore,
\begin{eqnarray*}
\int_0^T \bar{y}(t)^\top \bar{w}(t) dt & = & \int_0^T [b(t) - A(t)\bar{z}(t)]^\top \bar{w}(t) dt \\
& = & \int_0^T b(t)^\top \bar{w}(t) dt - \int_0^T \bar{z}(t)^\top A(t)^\top \bar{w}(t) dt \\
& = & \int_0^T b(t)^\top \bar{w}(t) dt - \int_0^T \bar{z}(t)^\top c(t) dt \\
& = & G(\bar{w}) - F(\bar{z}) = G(\tilde{w}) - F(\bar{z}) = 0.
\end{eqnarray*}
As $\bar{y}(t)^\top \bar{w}(t) \leq 0$ a.e. in $[0,T]$, (\ref{FC3}) holds.

Reciprocally, if there exists $\bar{y} \in L^{\infty}([0,T];\mathbb{R}^m)$ such that (\ref{FC1})-(\ref{FC3}) hold, it is clear that $\bar{z} \in \Omega$ and that $\bar{w} \in \Theta$. Moreover,
\begin{eqnarray*}
G(\bar{w}) & = & \int_0^T b(t)^\top \bar{w}(t) dt = \int_0^T [A(t)\bar{z}(t) + \bar{y}(t)]^\top \bar{w}(t) dt \\
& = & \int_0^T \bar{z}(t)^\top A(t)^\top \bar{w}(t) dt + \int_0^T \bar{y}(t)^\top \bar{w}(t) dt = \int_0^T \bar{z}(t)^\top c(t) dt \\
& = & F(\bar{z}).
\end{eqnarray*}
By Theorem \ref{DualidadeForte}, we conclude that $\bar{z} \in \Omega$ and that $\bar{w} \in \Theta$ are optimal solutions of (CLP) and (CDP).
\end{proof}

Note that the assumption that $\beta$-A is satisfied, in theorem above, is needed only to establish that conditions (\ref{FC1})-(\ref{FC3}) are necessary when we have optimal solutions of primal and dual problems. The other statement holds without such an assumption.

\begin{example}
We revisit Example \ref{Exemplo} one more time. Let us calculate the slack variables $\bar{y}_i(t) = b_i(t) - a_i(t)^\top \bar{z}(t)$, $i=1,2,3,4,5$. We have that
\begin{align*}
& \bar{y}_1(t) = \bar{z}_2(t) = \frac{1}{4}+\frac{5}{8}t ~ \mbox{a.e. in} ~ [0,2], \\
& \bar{y}_2(t) = \bar{z}_1(t) = \left\{ \begin{array}{l}
\frac{11}{4}-\frac{5}{8}t ~ \mbox{a.e. in} ~ [0,1], \\
\frac{1}{4}+\frac{5}{8}t ~ \mbox{a.e. in} ~ (1,2],
\end{array} \right. \\
& \bar{y}_3(t) = \left\{ \begin{array}{l} \bar{z}_1(t) - \bar{z}_2(t) = \frac{5}{2} - \frac{5}{4}t ~ \mbox{a.e. in} ~ [0,1], \\ -\bar{z}_1(t) + \bar{z}_2(t) = 0 ~ \mbox{a.e. in} ~ (1,2],  \end{array} \right. \\
& \bar{y}_4(t) = 3 - \bar{z}_1(t) - \bar{z}_2(t) = \left\{ \begin{array}{l} 0 ~ \mbox{a.e. in} ~ [0,1], \\ \frac{5}{2} - \frac{5}{4}t ~ \mbox{a.e. in} ~ (1,2],  \end{array} \right. \\
& \bar{y}_5(t) = \left( \frac{1}{4}+\frac{5}{8}t \right) - \bar{z}_2(t) = 0 ~ \mbox{a.e. in} ~ [0,2].
\end{align*}
We see that $\bar{z}$ and $\bar{y}$ satisfy (\ref{FC1}) by construction. Now we calculate $\bar{w}$ so that (\ref{FC3}) is valid. From the expressions for $\bar{y}_i$ above, (\ref{FC3}) is satisfied when
\begin{align*}
& \bar{w}_1(t) = \bar{w}_2(t) = 0 ~ \mbox{a.e. in} ~ [0,2], \\
& \bar{w}_3(t) = 0 ~ \mbox{a.e. in} ~ [0,1], \\
& \bar{w}_4(t) = 0 ~ \mbox{a.e. in} ~ (1,2]. 
\end{align*}
By using this information, we obtain from the constraints of the dual problem (see in Example \ref{Exemplo2}) that
\begin{align*}
& \bar{w}_3(t) = -t+1 ~ \mbox{a.e. in} ~ (1,2], \\
& \bar{w}_4(t) = t-1 ~ \mbox{a.e. in} ~ [0,1], \\
& \bar{w}_5(t) = -t ~ \mbox{a.e. in} ~ [0,2].
\end{align*}
Then, 
\begin{align*}
& \bar{w}_1(t) = \bar{w}_2(t) = 0 ~ \mbox{a.e. in} ~  [0,2], \\
& \bar{w}_3(t) = \left\{ \begin{array}{l} 0 ~ \mbox{a.e. in} ~ [0,1], \\ -t+1 ~ \mbox{a.e. in} ~ (1,2], \end{array} \right. \\
& \bar{w}_4(t) = \left\{ \begin{array}{l} t-1 ~ \mbox{a.e. in} ~ [0,1], \\ 0 ~ \mbox{a.e. in} ~ (1,2], \end{array} \right. \\
& \bar{w}_5(t) = -t ~ \mbox{a.e. in} ~ [0,2].
\end{align*}
Thus, by construction, (\ref{FC2}) is satisfied. Therefore, conditions (\ref{FC1})-(\ref{FC3}) are valid at $\bar{z}$, $\bar{y}$ and $\bar{w}$. It follows from Theorem \ref{FolgasCompl} that $\bar{w}$ is an optimal solution of the dual problem. 
\end{example}

\begin{theorem}
If the primal problem (CLP) has an optimal solution in which $\beta$-A is satisfied for some $\beta > 0$, then the dual problem (CDP) has an optimal solution.
\end{theorem}

\begin{proof}
If (CLP) has an optimal solution, we know from Corollaries \ref{CorolarioMultLag} and \ref{CorolarioMultLag2} that (CDP) also has an optimal solution.  
\end{proof}

\begin{theorem} \label{Teorema13}
Assume that hypothesis H is satisfied. If the dual problem (CDP) has an optimal solution, then the primal problem (CLP) has an optimal solution.
\end{theorem}

\begin{proof}
If (CDP) has an optimal solution, say $\bar{w}$, it follows from Theorem 4.1 in Monte and de Oliveira \cite{monte:2019}, that there exist $u \in L^\infty([0,T];\mathbb{R}^n)$ and $v \in L^\infty([0,T];\mathbb{R}^m)$ such that $b(t) + A(t)u(t) - v(t) = 0$, $v(t) \geq 0$ and $v(t)^\top \bar{w}(t) = 0$ a.e. in $[0,T]$ (hypothesis H implies that the full rank assumption in \cite{monte:2019} is satisfied). The pair $(\bar{z},\bar{y}) = (-u,v)$ satisfies (\ref{FC1}). It is clear that $\bar{w}$ satisfies (\ref{FC2}) and that $v(t)^\top \bar{w}(t) = 0$ a.e. in $[0,T]$ implies that (\ref{FC3}) also holds. It follows from Theorem \ref{FolgasCompl} (Complementary Slackness Theorem) that $\bar{z}$ is an optimal solution of (CLP). 
\end{proof}

\section{Conclusions}

Optimality conditions, as well as duality theory for continuous-time linear optimization problems with inequality constraints, was carried out. Classical duality results such as Karush-Kuhn-Tucker conditions, weak and strong duality properties, and the complementary slackness theorem were established, where regularity conditions were necessary. Such regularity conditions are less restrictive than those encountered in the literature. Thus, this work gives theoretical contributions to a class of continuous-time optimization problems. As already mentioned in the introduction, such problems were introduced in the fifties by Bellman in \cite{bellman:1953} as, according to Bellman himself, an interesting and significant class of production and allocation problems, called ``bottleneck problems''.

It is worth mentioning that a quick analysis will show that all proofs work in a more general framework of an abstract measure space. In particular, we can recover the known results from classical linear programming. 

The study of the continuous-time linear optimization problem in which the constraints are of the form
\begin{eqnarray*}
&& A(t)z(t) \leq b(t) + \int_0^t K(t,s)z(s)ds ~ \mbox{a.e. in} ~ [0,T], \\
&& z(t) \geq 0 ~ \mbox{a.e. in} ~ [0,T],
\end{eqnarray*}
is going to be a topic of future work. One other topic of future work is going to be deriving optimality conditions for the dual problem (CDP) under a more appropriate constraint qualification, so that Theorems \ref{Teorema9} and \ref{Teorema13} are valid down a less restrictive regularity condition other than hypothesis (H).

\bibliographystyle{plain}
\bibliography{LinearContTimeOptim}

\end{document}